\documentclass[12pt]{article}
\usepackage{amsmath}
\usepackage{amssymb}
\usepackage{vatola}

\def\q{\quad}
\def\qq{\qquad}
\def\qtq#1{\q\t{#1}\q}
\def\mod#1{\ (\text{\rm mod}\ #1)}
\def\t{\text}
\def\f{\frac}
\def\e{\equiv}
\def\b{\binom}
\def\sls#1#2{(\f{#1}{#2})}
 \def\ls#1#2{\big(\f{#1}{#2}\big)}
\def\Ls#1#2{\Big(\f{#1}{#2}\Big)}
\def\ap{\langle a\rangle_p}

\def\qp#1{q_p(#1)}

\let \pro=\proclaim
\let \endpro=\endproclaim

\begin{document}
 \centerline {\bf
Congruences for sums involving products of three binomial
coefficients }
\par\q\newline
\centerline{Zhi-Hong Sun}\newline
\centerline{School of Mathematics
and Statistics}
\centerline{Huaiyin Normal University}
\centerline{Huaian, Jiangsu 223300, P.R. China} \centerline{Email:
zhsun@hytc.edu.cn} \centerline{Homepage:
http://maths.hytc.edu.cn/szh1.htm}
 \abstract{ \par Let $p>3$ be a prime, and let $a$ be a rational $p$-adic integer,
  using WZ method we establish the congruences modulo
$p^3$ for
$$\sum_{k=0}^{p-1} \binom ak\binom{-1-a}k\binom{2k}k\frac {w(k)}{4^k},$$ where
$$w(k)=1,\frac 1{k+1},\frac 1{(k+1)^2},\frac 1{(k+1)^3},\frac 1{2k-1},
\frac 1{k+2},\frac 1{k+3},k,k^2,k^3,\frac 1{a+k},\frac 1{a+k-1}.$$
As consequences, taking $a=-\frac 12,-\frac 13,-\frac 14,-\frac 16$
we deduce many congruences modulo $p^3$ and so solve some
conjectures posed by the author earlier.
 \par\q
 }
 \endabstract

\section*{1. Introduction}
\par\q For $a\in\Bbb Z$ and given
odd prime $p$ let $\sls ap$ denote the Legendre symbol.  For
positive integers $a,b$ and $n$, if $n=ax^2+by^2$ for some integers
$x$ and $y$, we briefly write that $n=ax^2+by^2$. Let $p>3$ be a
prime. In 1987, Beukers[B] conjectured a congruence equivalent to
$$\sum_{k=0}^{p-1}\f{\b{2k}k^3}{64^k}\e
\cases 4x^2-2p\mod{p^2}&\t{if $p=x^2+4y^2\e 1\mod 4$,}
\\ 0\mod{p^2}&\t{if $p\e 3\mod 4$.}\endcases$$
This congruence was proved by several authors including
Ishikawa[I]($p\e 1\mod 4$), van Hamme[H]($p\e 3\mod 4$) and
Ahlgren[A]. Combining the results in [LR], [S10] and [T3], in [S14]
the author stated that
$$\sum_{k=0}^{p-1}\f{\b{2k}k^3}{64^k}\e
\cases 4x^2-2p-\f{p^2}{4x^2}\mod{p^3}\qq\qq\q\ \t{if $p=x^2+4y^2\e
1\mod 4$,}\\-\f {p^2}4\b{\f{p-3}2}{\f{p-3}4}^{-2}\e
-p^2\b{\f{p-1}2}{\f{p-3}4}^{-2}\mod{p^3}\qq\q\q\t{if $p\e 3\mod
4$.}\endcases\tag 1.1$$
\par Let $p>3$ be a prime. In 2003, Rodriguez-Villegas[RV] posed 22
conjectures on supercongruences modulo $p^2$. In particular, the
following congruences are equivalent to conjectures due to
Rodriguez-Villegas:
$$\align &\sum_{k=0}^{p-1}\f{\b{2k}k^2\b{3k}k}{108^k}\e
\cases 4x^2-2p\mod{p^2}&\t{if $p=x^2+3y^2\e 1\mod 3$,}
\\0\mod{p^2}&\t{if $p\e 2\mod 3$,}
\endcases
\\&\sum_{k=0}^{p-1}\f{\b{2k}k^2\b{4k}{2k}}{256^k}\e
\cases 4x^2-2p\mod{p^2}&\t{if $p=x^2+2y^2\e 1,3\mod 8$,}
\\0\mod{p^2}&\t{if $p\e 5,7\mod 8$,}
\endcases
\\&\Ls p3\sum_{k=0}^{p-1}\f{\b{2k}k\b{3k}k\b{6k}{3k}}{12^{3k}}\e
\cases 4x^2-2p\mod{p^2}&\t{if $p=x^2+4y^2\e 1\mod 4$,}
\\0\mod{p^2}&\t{if $p\e 3\mod 4$.}
\endcases\endalign$$
  These conjectures have been solved  by
Mortenson[M] and Zhi-Wei Sun[Su]. In 2018, J.C. Liu[Liu] conjectured
the congruences for
$$\sum_{k=0}^{p-1}\f{\b{2k}k^2\b{3k}k}{108^k},\q
\sum_{k=0}^{p-1}\f{\b{2k}k^2\b{4k}{2k}}{256^k},\q
\sum_{k=0}^{p-1}\f{\b{2k}k\b{3k}k\b{6k}{3k}}{1728^k}\mod {p^3}$$ in
terms of $p-$adic Gamma functions. In [S11], the author conjectured
that
$$\align &\qq\sum_{k=0}^{p-1}\f{\b{2k}k^2\b{3k}k}{108^k}\e
\cases 4x^2-2p-\f{p^2}{4x^2}\mod{p^3}&\t{if $p=x^2+3y^2\e 1\mod 3$,}
\\-\f{p^2}2\b{(p-1)/2}{(p-5)/6}^{-2}\mod{p^3}&\t{if $p\e 2\mod 3$,}
\endcases\tag 1.2
\\&\qq\sum_{k=0}^{p-1}\f{\b{2k}k^2\b{4k}{2k}}{256^k}\e
\cases 4x^2-2p-\f{p^2}{4x^2}\mod{p^3}&\t{if $p=x^2+2y^2\e 1,3\mod
8$,}
\\-\f{p^2}3\b{(p-1)/2}{[p/8]}^{-2}\mod{p^3}&\t{if $p\e 5,7\mod 8$,}
\endcases\tag 1.3
\\&\qq\Ls p3\sum_{k=0}^{p-1}\f{\b{2k}k\b{3k}k\b{6k}{3k}}{12^{3k}}\e
\cases 4x^2-2p-\f{p^2}{4x^2}\mod{p^3}\q \t{if $p=x^2+4y^2\e 1\mod
4$,}
\\\f 5{12}p^2\b{(p-3)/2}{(p-3)/4}^{-2}\mod{p^3}\qq\qq\t{if $p\e 3\mod 4$,}
\endcases\tag 1.4\endalign$$
where $[a]$ is the greatest integer not exceeding $a$.  It is easy
to see that (see [S4-S7])
$$\align &\b{-\f 12}k=\f{\b{2k}k}{(-4)^k},\q\b {-\f 13}k\b{-\f
23}k=\f{\b{2k}k\b{3k}k}{27^k},
\\&\b {-\f 14}k\b{-\f 34}k=\f{\b{2k}k\b{4k}{2k}}{64^k},
\q\b {-\f 16}k\b{-\f 56}k=\f{\b{3k}{k}\b{6k}{3k}}{432^k}.
\endalign$$
Thus, a natural and general problem is to determine
$\sum_{k=0}^{p-1} \b ak\b{-1-a}k\b{2k}k\f 1{4^k}$ modulo $p^3$,
where $p>3$ is a prime and $a$ is a rational $p$-adic integer.
\par For a prime $p$ let $\Bbb Z_p$ be the set of rational numbers whose denominator is not divisible by
$p$. For $a\in \Bbb Z_p$ let $\ap$ be determined by
$\ap\in\{0,1,\ldots,p-1\}$ and $a\e \ap\mod p$. In [S7], the author
showed that for any given odd prime $p$ and $a\in\Bbb Z_p$,
$$\sum_{k=0}^{p-1} \b ak\b{-1-a}k\b{2k}k\f 1{4^k}\e 0\mod {p^2}
\qtq{for}\ap\e 1\mod 2.$$ For an odd prime $p$ and $x\in\Bbb Z_p$,
 the $p$-adic Gamma function $\Gamma_p(x)$ is
  defined by
  $$\Gamma_p(0)=1,\q \Gamma_p(n)=(-1)^n\prod_{
  \substack{k\in\{1,2,\ldots,n-1\}\\p\nmid
  k}}k\qtq{for}n=1,2,3,\ldots$$
  and
  $$\Gamma_p(x)=\lim_{\substack{n\in\{0,1,\ldots\}\\ |x-n|_p\rightarrow
  0}}
  \Gamma_p(n).$$
  In [PTW], Pan, Tauraso and Wang established a result equivalent to
  $$\sum_{k=0}^{p-1} \b ak\b{-1-a}k\b{2k}k\f 1{4^k}\e
  \cases \f{\Gamma_p\sls
  12^2}{\Gamma_p\sls{2+a}2^2\Gamma_p(\f{1-a}2)^2}\mod {p^3}&\t{if
  $2\mid \ap$,}
  \\\f{a'(a'+1)}4p^2\f{\Gamma_p\sls
  12^2}{\Gamma_p\sls{2+a}2^2\Gamma_p(\f{1-a}2)^2}\mod {p^3}&\t{if
  $2\nmid \ap$,}\endcases\tag 1.5$$
  where $a'=(a-\ap)/p$.
But they only gave the proof for the case $2\mid \ap$, and their
method is somewhat complicated. In the case $2\mid \ap$, (1.5) was
first conjectured by Liu in [Liu]. Using (1.5) and Jacobi sums,
recently Mao [Mao] proved (1.2) and (1.3). We also note that Guo
[Guo] established three congruences modulo $p^3$ concerning the sums
in (1.2)-(1.4) via $q$-congruences. For example, Guo showed that for
any prime $p\e 5\mod 6$,
$$\sum_{k=0}^{p-1}\f{\b{2k}k^2\b{3k}k}{108^k}
\e \f{2p\b{-2/3}{(2p-1)/3}}{3\b{-7/6}{(2p-1)/3}} +
\f{p\b{-5/6}{(p-2)/3}}{3\b{-4/3}{(p-2)/3}}\mod {p^3}.$$ In [Su],
Z.W. Sun showed that
$$\align &\sum_{k=0}^{p-1}\f{\b{2k}k^3}{64^k(k+1)}
\e 4x^2-2p\mod {p^2}\qtq{for}p=x^2+4y^2\e 1\mod 4,
\\&\sum_{k=0}^{p-1}\f{\b{2k}k^2\b{3k}{k}}{108^k(k+1)}
\e 4x^2-2p\mod {p^2}\qtq{for}p=x^2+3y^2\e 1\mod 3,
\\&\sum_{k=0}^{p-1}\f{\b{2k}k^2\b{4k}{2k}}{256^k(k+1)}
\e 4x^2-2p\mod {p^2}\qtq{for}p=x^2+2y^2\e 1,3\mod 8,
\\&\Ls p3\sum_{k=0}^{p-1}\f{\b{2k}k\b{3k}k\b{6k}{3k}}{1728^k(k+1)} \e 4x^2-2p\mod
{p^2}\qtq{for}p=x^2+4y^2\e 1\mod 4.
\endalign$$
In [T3], Tauraso showed that
$$\align&\sum_{k=0}^{(p-1)/2}\f{\b{2k}k^3}{64^k(k+1)^3}
\\&\e\cases 8-\f{24p^2}{\Gamma_p\sls 14^4}\e 8+\f{6p^2}{x^2}\mod {p^3}
\qq\qq\qq\qq\t{if $p=x^2+4y^2\e 1\mod 4$,}
\\8-\f{384}{\Gamma_p\sls 14^4}\e
8-\f{96}{2^{p-1}(1+2p+(3-\f
12E_{p-3})p^2)}\b{\f{p-3}2}{\f{p-3}4}^2\mod {p^3} \q\t{if $4\mid
p-3$,}
\endcases\endalign$$
where  $\{E_n\}$ are the Euler numbers given by
$$E_{2n-1}=0,\q E_0=1,\q E_{2n}=-\sum_{k=1}^n\b{2n}{2k}E_{2n-2k}\
(n=1,2,3,\ldots).$$ In [S11], [S13] and [S14], the author posed
numerous conjectures on congruences modulo $p^3$ for the sums
$$\sum_{k=0}^{p-1}\f{w(k)\b{2k}k^3}{m^k},\q
\sum_{k=0}^{p-1}\f{w(k)\b{2k}k^2\b{3k}k}{m^k},\q
\sum_{k=0}^{p-1}\f{w(k)\b{2k}k^2\b{4k}{2k}}{m^k},\q
\sum_{k=0}^{p-1}\f{w(k)\b{2k}k\b{3k}k\b{6k}{3k}}{m^k},$$ where $m$
is an integer not divisible by $p$ and $w(k)\in\{1,\f 1{k+1},\f
1{2k-1},\f 1{(k+1)^2},k,k^2,k^3\}$.

\par Let $p>3$ be a prime and $a\in\Bbb Z_p$. Inspired by the above work, in
this paper, using WZ method we establish the congruences modulo
$p^3$ for
$$\sum_{k=0}^{p-1} \b ak\b{-1-a}k\b{2k}k\f {w(k)}{4^k},$$ where
$$w(k)\in\Big\{1,\f 1{k+1},\f 1{(k+1)^2},\f 1{(k+1)^3},\f 1{2k-1},\f
1{k+2},\f 1{k+3}, k,k^2,k^3,\f 1{a+k},\f 1{a+k-1}\Big\}.$$ Our
approach is natural and elementary. As consequences, taking $a=-\f
12,-\f 13,-\f 14,-\f 16$ we deduce many congruences modulo $p^3$ and
so solve some conjectures in [S8,S10,S11]. For instance, for $p\e
1\mod 4$ and so $p=x^2+4y^2$ we prove that
$$\align&\sum_{k=0}^{(p-1)/2}\f{\b{2k}k^3}{64^k(k+1)}
\e 4x^2-2p\mod {p^3},
\\&\sum_{k=0}^{(p-1)/2}\f{\b{2k}k^3}{64^k(k+2)}\e
 \f{52}{27}x^2-\f{26}{27}p-\f{p^2}{27x^2}\mod {p^3},
\\&\sum_{k=0}^{p-1}\f{\b{2k}k^3}{64^k(2k-1)}
\e -2x^2+p+\f{p^2}{4x^2} \mod {p^3},
\\&\sum_{k=0}^{p-1}\f{\b{2k}k^3}{64^k(2k-1)^2}
\e 2x^2-p-\f{p^2}{2x^2} \mod {p^3},
\endalign$$
for $p\e 1\mod 3$, $p>7$ and so $p=x^2+3y^2$ we show that
$$\align
&\sum_{k=0}^{p-2}\f{\b{2k}k^2\b{3k}{k}}{108^k(k+1)}\e 4x^2-2p\mod
{p^3}, \\&\sum_{k=0}^{p-2}\f{\b{2k}k^2\b{3k}k}{108^k(k+1)^2} \e
8x^2-4p+\f{9p^2}{8x^2} \mod {p^3}
\\&\sum_{k=0}^{p-2}\f{\b{2k}k^2\b{3k}k}{108^k(k+1)^3} \e
9-2x^2+p+\f{117p^2}{16x^2} \mod {p^3},
\\&\sum_{k=0}^{p-3}\f{\b{2k}k^2\b{3k}k}{108^k(k+2)} \e
\f{29}{15}x^2-\f{29}{30}p-\f{3p^2}{80x^2} \mod {p^3},
\\&\sum_{k=0}^{p-4}\f{\b{2k}k^2\b{3k}k}{108^k(k+3)} \e
\f{32}{25}x^2-\f{16}{25}p-\f{99p^2}{2800x^2} \mod {p^3},
\\&\sum_{k=0}^{p-1}\f{\b{2k}k^2\b{3k}k}{108^k(2k-1)}\e
-\f{20}9x^2+\f{10}9p+\f {p^2}{4x^2}\mod {p^3},
\\&\sum_{k=0}^{p-1}\f{k\b{2k}k^2\b{3k}{k}}{108^k}
\e -\f 89x^2+\f 49p+\f{p^2}{9x^2}\mod {p^3},
\\&\sum_{k=0}^{p-1}\f{k^2\b{2k}k^2\b{3k}{k}}{108^k}
\e \f {32}{243}x^2-\f {16}{243}p-\f{17p^2}{486x^2}\mod {p^3},
\\&\sum_{k=0}^{p-1}\f{k^3\b{2k}k^2\b{3k}{k}}{108^k}
\e  \f 1{10935}\Big(16x^2-8p+\f{113p^2}{2x^2}\Big)\mod {p^3}.
\endalign$$
We also pose some challenging conjectures on the congruences modulo
$p^3$ for the sum $\sum_{k=0}^{(p-1)/2}w(k)\f{\b{2k}k^3}{m^k}$,
where $m\in\{1,-8,16,-64,256,-512,4096\}$ and $w(k)\in\{\f 1{k+2},\f
1{k+3},$ $\f 1{(k+1)^2},\f 1{(k+2)^2},\f 1{(k+1)^3}\}$.

 \par In addition to the above notation,
throughout this paper we use the following notations.
   Let $H_0=H_0^{(2)}=0$. For $n\ge 1$ let $H_n=1+\f 12+\cdots+\f 1n$ and
   $H_n^{(2)}=1+\f 1{2^2}+\cdots+\f 1{n^2}$. For an odd prime and
   $a\in\Bbb Z_p$ set $q_p(a)=(a^{p-1}-1)/p$ and
   $$\align &R_1(p)=(2p+2-2^{p-1})\b{(p-1)/2}{[p/4]}^2,
   \\&R_2(p)=(5-4(-1)^{\f{p-1}2})\Big(1+(4+2(-1)^{\f{p-1}2})p
-4(2^{p-1}-1)-\f p2\sum_{k=1}^{[p/8]} \f 1k\Big)\b{\f{p-1}2}{[\f
p8]}^2,
\\&R_3(p)=\Big(1+2p+\f
43(2^{p-1}-1)-\f 32(3^{p-1}-1)\Big) \b{(p-1)/2}{[p/6]}^2,
\\&R_7(p)=\sum_{k=0}^{(p-1)/2}\f{\b{2k}k^3}{k+1}.
\endalign$$
Let
$$S_n(a)=\sum_{k=0}^{n-1}\b ak\b {-1-a}k\b{2k}k\f 1{4^k}\q(n=1,2,3,\ldots),$$
 and let $\{U_n\}$ be the sequence given by
$$U_{2n-1}=0,\q U_0=1,\q U_{2n}=-2\sum_{k=1}^n\b{2n}{2k}U_{2n-2k}\
(n=1,2,3,\ldots).$$

\section*{2. The congruence for
$\sum_{k=0}^{p-1}\b{a}k\b{-1-a}k\b{2k}{k}\f 1{4^k}$ modulo $p^3$}
\par For $k=0,1,2,\ldots$ set
$$F(a,k)=\b ak\b{-1-a}k\b{2k}{k}\f 1{4^k}$$
and
$$G(a,k)=(a+2)(2a+3)\f{k}{4^{k-1}(a+1+k)}\b{2k-1}{k-1}\b{a+1}{k-1}
\b{-3-a}{k-1}.$$ It is easy to check that
$$(a+2)^2F(a+2,k)-(a+1)^2F(a,k)=G(a,k+1)-G(a,k).$$
Thus,
$$\align &(a+2)^2S_n(a+2)-(a+1)^2S_n(a)
\\&=(a+2)^2\sum_{k=0}^{n-1}F(a+2,k)-(a+1)^2\sum_{k=0}^nF(a,k)
\\&=\sum_{k=0}^{n-1}(G(a,k+1)-G(a,k))=G(a,n)-G(a,0)=G(a,n).\endalign$$
That is,
$$\aligned &(a+2)^2S_n(a+2)-(a+1)^2S_n(a)
\\&=(a+2)(2a+3)\f{n}{4^{n-1}(a+1+n)}\b{2n-1}{n-1}\b{a+1}{n-1}\b{-3-a}{n-1}.
\endaligned\tag 2.1$$

\pro{Lemma 2.1} Let $p>3$ be a prime, $a\in\Bbb Z_p$, $a\not\e
-1\mod p$ and $a'=(a-\ap)/p$. Then
$$\align &(a+2)^2S_{p}(a+2)-(a+1)^2S_{p}(a)\\&\e \cases \big(\f 1{a+1}+\f
1{a+2}\big)a'(a'+1)p^3\mod {p^4}&\t{if $\ap<p-2$,}
\\(a+2)p\mod{p^3}&\t{if $\ap=p-2$.}\endcases\endalign$$
\endpro
Proof. From (2.1) we see that
$$\align &(a+2)^2S_p(a+2)-(a+1)^2S_p(a)
\\&=(a+2)(2a+3)\f{p}{4^{p-1}(a+1+p)}\b{2p-1}{p}\b{a+2-1}{p-1}\b{-a-2-1}{p-1}.
\endalign$$
We first assume $\ap<p-2$. Then $a+2-\langle
a+2\rangle_p=a-\ap=a'p$. From [S9, Lemma 2.2] we know that
$$\b{a+2-1}{p-1}\b{-a-2-1}{p-1}\e \f{a'(a'+1)p^2}{(a+2)^2}\mod
{p^3}.\tag 2.2$$ Hence,
$$\align &(a+2)^2S_{p}(a+2)-(a+1)^2S_{p}(a)\\&\e
(a+2)(2a+3)\f p{4^{p-1}(a+1+p)}\b{2p-1}{p-1}\f{a'(a'+1)p^2}{(a+2)^2}
\\&\e \Big(\f 1{a+1}+\f
1{a+2}\Big)a'(a'+1)p^3\mod {p^4}.\endalign$$
\par Now we assume that $\ap=p-2$. Then $a+2=p(a'+1)$. Appealing to [S12, (2.6)],
$$\align \b{a+2-1}{p-1}\b{-a-2-1}{p-1}
&=\b{p-1+a'p}{p-1}\b{p-1-(a'+2)p}{p-1}\\& \e
(1+a'pH_{p-1})(1-(a'+2)pH_{p-1})\e 1\mod {p^2}.\endalign$$ Hence,
$$\align &(a+2)^2S_{p}(a+2)-(a+1)^2S_{p}(a)\\&\e
(a+2)(2a+3)\f p{4^{p-1}(a+1+p)}\b{2p-1}{p-1}\e (a+2)p\mod
{p^3}.\endalign$$  This completes the proof.

\pro{Lemma 2.2} Let $p$ be an odd prime and $t\in\Bbb Z_p$. Then
$$\sum_{k=0}^{p-1}\b{pt}k\b{-1-pt}k\b{2k}k\f 1{4^k}
\e 1-2t(2^{p-1}-1)+(2t^2+t)(2^{p-1}-1)^2\mod {p^3}.$$
\endpro
Proof. Clearly,
$$\align &\sum_{k=0}^{p-1}\b{pt}k\b{-1-pt}k\b{2k}k\f 1{4^k}
\\&=1+\sum_{k=1}^{p-1}\f{pt}{pt-k}\cdot\f{((-1)^2-p^2t^2)((-2)^2-p^2t^2)
\cdots((-k)^2-p^2t^2)}{k!^2}\cdot\f{\b{2k}k}{4^k}
\\&\e 1-\sum_{k=1}^{p-1}\f{pt(pt+k)}{k^2}\cdot\f{\b{2k}k}{4^k}
\mod {p^3}.\endalign$$ By [S8, Remark 3.1] or [T1],
$$\sum_{k=1}^{p-1}\f{\b{2k}k}{k\cdot 4^k}
=\sum_{k=1}^{p-1}\f{(-1)^k}k\b{-1/2}k \e 2q_p(2)-pq_p(2)^2\mod
{p^2}.\tag 2.3$$ Taking $x=\f 14$ in [T2, (9)] and then applying
[S2, Theorem 4.1] we see that
$$\sum_{k=1}^{p-1}\f{\b{2k}k}{k^2\cdot 4^k}\e 4\sum_{k=1}^{p-1}\f
1{k^2\cdot 2^k}\e -2q_p(2)^2\mod p.\tag 2.4$$ It then follows that
$$\sum_{k=0}^{p-1}\b{pt}k\b{-1-pt}k\b{2k}k\f 1{4^k}
\e 1-pt(2q_p(2)-pq_p(2)^2)-p^2t^2(-2q_p(2)^2)\mod {p^3},$$ which
yields the result.

\pro{Lemma 2.3} Let $p$ be an odd prime,
$n\in\{1,2,\ldots,\f{p-1}2\}$ and $t\in\Bbb Z_p$. Then
$$\align &\b{\f{p-1}2+pt}n\\&\e \b{\f{p-1}2}n\Big(1-2pt\sum_{k=1}^n\f 1{2k-1}
+2p^2t\Big(t\Big(\sum_{k=1}^n\f 1{2k-1}\Big)^2-(t+1)\sum_{k=1}^n\f
1{(2k-1)^2}\Big)\mod {p^3}.\endalign$$
\endpro

Proof. For $m\in\{n+1,\ldots,p-1\}$ we see that
$$\aligned &\b{m+pt}n=\f{(pt+m)(pt+m-1)\cdots(pt+m-(n-1))}{n!}\\&\e
\b mn\Big(1+pt\sum_{m-n+1\le k\le m}\f 1k+p^2t^2\sum_{m-n+1\le
i<j\le m}\f 1{ij}\Big)
\\&=\b mn\Big(1+pt(H_m-H_{m-n})+\f
12p^2t^2\big((H_m-H_{m-n})^2-(H_m^{(2)}-H_{m-n}^{(2)})\big)\Big)
\mod {p^3}.
\endaligned\tag 2.5$$
For given positive integer $r$ we have
$$\align &\sum_{k=1}^{(p-1)/2}\f 1{k^r}-\sum_{k=1}^{(p-1)/2-n}\f
1{k^r} \\&=\sum_{k=1}^n\f 1{\ls{p-(2k-1)}2^r}=2^r\sum_{k=1}^n
\f{(p+2k-1)^r}{(p^2-(2k-1)^2)^r} \e (-2)^r\sum_{k=1}^n
\f{(2k-1)^r+rp(2k-1)^{r-1}}{(2k-1)^{2r}} \\&=(-2)^r\sum_{k=1}^n\f
1{(2k-1)^r}+(-2)^rrp\sum_{k=1}^n\f 1{(2k-1)^{r+1}}\mod {p^2}.
\endalign$$
Hence,
$$\align &H_{\f{p-1}2}-H_{\f{p-1}2-n}\e -2\sum_{k=1}^n\f
1{2k-1}-2p\sum_{k=1}^n\f 1{(2k-1)^2}\mod {p^2},
\\&\sum_{k=1}^{(p-1)/2}\f 1{k^2}-\sum_{k=1}^{(p-1)/2-n}\f 1{k^2}
\e 4\sum_{k=1}^n\f 1{(2k-1)^2}\mod p.\endalign$$ Now, from the above
we deduce that
$$\align \b{\f{p-1}2+pt}n&\e \b{\f{p-1}2}n\Big(1+pt\Big(-2
\sum_{k=1}^n\f 1{2k-1}-2p\sum_{k=1}^n\f 1{(2k-1)^2}\Big)\\&\q+\f
12p^2t^2 \Big(\Big(-2\sum_{k=1}^n\f 1{2k-1}\Big)^2-4\sum_{k=1}^n\f
1{(2k-1)^2}\Big)\mod {p^3},\endalign$$ which yields the result.

\pro{Theorem 2.1} Let $p$ be an odd prime, $a\in\Bbb Z_p$, $a\not\e
0,-1\mod p$ and $a'=(a-\ap)/p$. If $2\mid \ap$ and $\ap=2n$, then
$$\align S_p(a)&\e \f{\b{(a-1)/2}{n}^2}{\b{a/2}{n}^2}
\big(1-2a'(2^{p-1}-1)+a'(2a'+1)(2^{p-1}-1)^2)
\\&\e
\b{(p-1)/2}n^2\Big(1+p\big((2a'+2)H_{2n}-(2a'+1)H_n-2a'\qp2\big)
\\&\q+\f{p^2}2 \Big(2a'\qp 2^2+\big((2a'+2)H_{2n}-(2a'+1)H_n-2a'\qp2\big)^2
\\&\q+\f 12(2{a'}^2-1)H_n^{(2)}+2(1-{a'}^2)H_{2n}^{(2)}\Big)\Big)\mod{p^3}.
\endalign$$ If $2\nmid \ap$,
then
$$S_p(a)\e\f{4^{\ap-1}\cdot a'(a'+1)p^2}{a^2\b{\ap-1}{\f{\ap-1}2}^2} \e
\f{a'(a'+1)p^2}{a^2\b{(p-1)/2}{(\ap-1)/2}^2} \mod {p^3}.$$
\endpro
Proof. Set $a'=(a-\ap)/p$. We first assume that $2\mid \ap$ and
$\ap=2n$. From Lemma 2.1 we see that
$$\align S_p(a)
&\e \f{(a-1)^2}{a^2}S_p(a-2)\e  \f{(a-1)^2}{a^2}\cdot
\f{(a-3)^2}{(a-2)^2}S_p(a-4)\e\cdots
\\&\e
\f{(a-1)^2(a-3)^2\cdots(a-2n+1)^2}{a^2(a-2)^2\cdots(a-2n+2)^2}S_p(a-2n)
\\&=\prod_{k=1}^{n}\f{(\f{a+1}2-k)^2}{(\f{a+2}2-k)^2}\cdot S_p(a-2n)
=\f{\b{(a-1)/2}{n}^2}{\b{a/2}{n}^2} S_p(a-2n)\mod {p^3}.
\endalign$$
By Lemma 2.2,
$$S_p(a-2n)=S_p(a'p)\e 1-2a'(2^{p-1}-1)+a'(2a'+1)(2^{p-1}-1)^2
\mod {p^3}.$$ Thus,
$$S_p(a)\e \b{(a-1)/2}{n}^2\b{a/2}{n}^{-2}\big(1-2a'p\qp 2+
a'(2a'+1)p^2\qp 2^2\big)\mod {p^3}.$$ Since
$(1+bp+cp^2)(1-bp+(b^2-c)p^2)\e 1\mod {p^3}$, appealing to (2.5) we
get
$$\align \b{a/2}n^{-1}\e \b{n+a'p/2}n^{-1}
&\e \Big(1+\f 12a'pH_n+\f 18{a'}^2p^2(H_n^2-H_n^{(2)})\Big)^{-1}
\\&\e 1-\f 12a'pH_n+\f 18{a'}^2p^2(H_n^2+H_n^{(2)})\mod {p^3}.
\endalign$$
By Lemma 2.3,
$$\align \b{\f{a-1}2}n&=(-1)^n\b{-\f{a-1}2+n-1}n=(-1)^n\b{\f{p-1}2-\f{a'+1}2p}n
\\&\e (-1)^n\b{\f{p-1}2}n\Big(1+(a'+1)p\sum_{k=1}^n\f 1{2k-1}
\\&\q+(a'+1)p^2\Big(\f{a'+1}2\Big(\sum_{k=1}^n\f
1{2k-1}\Big)^2+\f{1-a'}2\sum_{k=1}^n\f 1{(2k-1)^2}\Big)\mod {p^3}.
\endalign$$
Therefore,
$$ \align &(-1)^n\b{\f{a-1}2}n\b{\f a2}n^{-1}\b{\f{p-1}2}n^{-1}
\\&\e \Big(1+(a'+1)p\sum_{k=1}^n\f 1{2k-1}
+(a'+1)p^2\Big(\f{a'+1}2\Big(\sum_{k=1}^n\f
1{2k-1}\Big)^2+\f{1-a'}2\sum_{k=1}^n\f 1{(2k-1)^2}\Big)
\\&\q\times\Big(1-\f 12a'pH_n+\f 18{a'}^2p^2(H_n^2+H_n^{(2)})\Big)
\\&\e 1+p\Big((a'+1)\sum_{k=1}^n\f 1{2k-1}-\f {a'}2H_n\Big)
+p^2\Big(\f 18{a'}^2H_n^2+\f 18{a'}^2H_n^{(2)}
\\&\q+
\f{(a'+1)^2}2\Big(\sum_{k=1}^n\f 1{2k-1}\Big)^2
+\f{1-{a'}^2}2\sum_{k=1}^n\f 1{(2k-1)^2}-\f
12a'(a'+1)H_n\sum_{k=1}^n\f 1{2k-1}\Big)
\\&=1+p\Big((a'+1)\sum_{k=1}^n\f 1{2k-1}-\f {a'}2H_n\Big)
\\&\q+p^2\Big(\f 18\Big(a'H_n-2(a'+1)\sum_{k=1}^n\f 1{2k-1}\Big)^2
+\f 18{a'}^2H_n^{(2)}+\f{1-{a'}^2}2\sum_{k=1}^n\f 1{(2k-1)^2}\Big)
\mod {p^3}.\endalign$$ Note that $(1+bp+cp^2)^2\e
1+2bp+(b^2+2c)p^2\mod {p^3}$. From the above we derive that
$$ \align &\b{\f{a-1}2}n^2\b{\f a2}n^{-2}\b{\f{p-1}2}n^{-2}
\\&\e 1+p\Big(2(a'+1)\sum_{k=1}^n\f 1{2k-1}-a'H_n\Big)
+p^2\Big(2\Big((a'+1)\sum_{k=1}^n\f 1{2k-1}-\f {a'}2H_n\Big)^2
\\&\q+\f 14{a'}^2H_n^{(2)}+(1-{a'}^2)\sum_{k=1}^n\f 1{(2k-1)^2}\Big)
\\&=1+p\Big(2(a'+1)(H_{2n}-\f 12H_n)-a'H_n\Big)
+p^2\Big(2\Big((a'+1)(H_{2n}-\f 12H_n)-\f {a'}2H_n\Big)^2
\\&\q+\f 14{a'}^2H_n^{(2)}+(1-{a'}^2)(H_{2n}^{(2)}-\f 14H_n^{(2)}\Big)
\\&=1+p((2a'+2)H_{2n}-(2a'+1)H_n)+p^2\Big(\f
12((2a'+2)H_{2n}-(2a'+1)H_n)^2\\&\q+ (1-{a'}^2)H_{2n}^{(2)}
+\f{2{a'}^2-1}4H_n^{(2)}\Big) \mod {p^3}\endalign$$ and so
$$\align S_p(a)&\e\b{(p-1)/2}n^2(1-2a'p\qp
2+a'(2a'+1)p^2\qp 2^2)\\&\q\times\Big(1+p((2a'+2)H_{2n}-(2a'+1)H_n)
+p^2\Big(\f 12((2a'+2)H_{2n}-(2a'+1)H_n)^2\\&\q+
(1-{a'}^2)H_{2n}^{(2)} +\f{2{a'}^2-1}4H_n^{(2)}\Big)\Big)
\\&\e \b{(p-1)/2}n^2\Big(1+p(-2a'\qp 2+2(a'+1)H_{2n}-(2a'+1)H_n)+p^2
\Big(a'\qp 2^2\\&\q+\f 12\big((2a'+2)H_{2n}-(2a'+1)H_n-2a'\qp
2\big)^2\\&\q + (1-{a'}^2)H_{2n}^{(2)}
+\f{2{a'}^2-1}4H_n^{(2)}\Big)\Big) \mod {p^3},
\endalign$$
which yields the result in the case $2\mid \ap$.
\par Now assume that $2\nmid \ap$. By Lemma 2.1, for $\ap\le p-4$,
$$\align S_p(a)&\e \f{(a+2)^2}{(a+1)^2}S_p(a+2)\e
\f{(a+2)^2}{(a+1)^2}\cdot \f{(a+4)^2}{(a+3)^2}S_p(a+4)\e\cdots
\\&\e\prod_{k=0}^{\f{p-\ap}2-2}\f{(a+2k+2)^2}{(a+2k+1)^2}\cdot
S_p(p+a-\ap-2)
\\&=\prod_{k=0}^{\f{p-\ap}2-2}\f{(a+2k+1)^2(a+2k+2)^2}{(a+2k+1)^4}
\cdot S_p(p+a-\ap-2)
\\&=\f{(a+1)^2(a+2)^2\cdots(p+a-\ap-3)^2(p+a-\ap-2)^2}
{2^{4(\f{p-\ap}2-1)}(\f{a+1}2(\f{a+1}2+1)\cdots(\f{a+1}2+\f{p-\ap}2-2))^4}
\\&\q\times S_p(p+a-\ap-2) \mod {p^3}.
\endalign$$
By Lemmas 2.1 and 2.2,
$$\align((a'+1)p-1)^2 S_p((a'+1)p-2)&\e
(a'+1)^2p^2S_p((a'+1)p)-(a'+1)p^2 \\&\e
(a'+1)^2p^2-(a'+1)p^2=a'(a'+1)p^2 \mod {p^3}.\endalign$$ Hence,
$$S_p(p+a-\ap-2)=S_p((a'+1)p-2)\e a'(a'+1)p^2 \mod {p^3}.$$
Now, from the above we deduce that
$$\align S_p(a)&\e a'(a'+1)p^2\cdot\f{(\ap+1)^2(\ap+2)^2
\cdots(p-3)^2(p-2)^2}
{2^{4(\f{p-\ap}2-1)}(\f{\ap+1}2(\f{\ap+1}2+1)\cdots(\f{p-3}2))^4}
\\&=
a'(a'+1)p^2\cdot\f{(p-2)!^2\cdot\f{\ap-1}2!^4}{\ap!^2\cdot
2^{2(p-\ap)-4}\cdot \f{p-3}2!^4}
\\&=a'(a'+1)p^2\cdot\f{{\sls{p-1}2}^4}{(p-1)^2}
\b{p-1}{\f{p-1}2}^2\cdot\f
1{2^{2(p-\ap)-4}\cdot\sls{\ap+1}2^2\b{\ap}{\f{\ap-1}2}^2}
\\&\e\f{a'(a'+1)p^2}{2^{2(p-\ap)}\sls{a+1}2^2\b{\ap}{\f{\ap-1}2}^2}
\e \f{4^{\ap}\cdot a'(a'+1)p^2}{(a+1)^2\b{\ap}{\f{\ap-1}2}^2} \\&\e
\f{4^{\ap-1}\cdot a'(a'+1)p^2}{a^2\b{\ap-1}{\f{\ap-1}2}^2}
=\f{a'(a'+1)p^2}{a^2\b{-1/2}{(\ap-1)/2}^2} \e
\f{a'(a'+1)p^2}{a^2\b{(p-1)/2}{(\ap-1)/2}^2}
 \mod {p^3}.\endalign$$
For $\ap=p-2$, using Lemmas 2.1 and 2.2 we see that
$$\align S_p(a)&\e \f 1{(a+1)^2}\big((a+2)^2S_p(a+2)-(a+2)p\big)\e
(a+2)^2-(a+2)p\\&\e  \f{4^{\ap-1}\cdot
a'(a'+1)p^2}{a^2\b{\ap-1}{\f{\ap-1}2}^2}\e
\f{a'(a'+1)p^2}{a^2\b{(p-1)/2}{(\ap-1)/2}^2}\mod {p^3}.\endalign$$
This completes the proof.

\pro{Corollary 2.1} Let $p$ be a prime with $p>3$ and $p\e 3\mod 4$.
Then
$$\sum_{k=0}^{p-1}\f{\b{2k}k\b{3k}k\b{6k}{3k}}{1728^k}
\e\cases -\f{5p^2}{\b{(p-1)/2}{[p/12]}^2}\mod {p^3}&\t{if $p\e
7\mod{12}$,}
\\-\f{p^2}{5\b{(p-1)/2}{[p/12]}^2}\mod {p^3}&\t{if $p\e
11\mod{12}$.}\endcases$$
\endpro
Proof. Observe that $\b{-\f 16}k\b{-\f 56}k=\b{3k}k\b{6k}{3k}\f
1{432^k}$. Let $a=-\f 16$ or $-\f 56$ according as $p\e 7\mod {12}$
or $p\e 11\mod {12}$, and  $a'=(a-\ap)/p$. Then $\ap=[\f p6]$,
$\f{\ap-1}2=[\f p{12}]$ and $a'(a'+1)=-\f 5{36}$. Now, the result
follows from Theorem 2.1 immediately.

\pro{Lemma 2.4 (See [S1, Theorem 5.2] and [S2, Corollaries 3.3 and
3.7])} Let $p>3$ be a prime. Then
$$\align &H_{\f{p-1}2}\e -2\qp 2+p\qp 2^2\mod {p^2},
\q H_{\f{p-1}2}^{(2)}\e 0\mod p,
\\&H_{[\f p4]}\e -3\qp 2+\f 32p\qp 2^2-(-1)^{\f{p-1}2}pE_{p-3}\mod
{p^2},
\\&H_{[\f p4]}^{(2)}\e 4(-1)^{\f{p-1}2}E_{p-3}\mod p.
\endalign$$
\endpro
\par We remark that putting $a=-\f 12$ in Theorem 2.1 and then
applying Lemma 2.4 gives a natural and elementary proof of (1.1).

 \pro{Lemma 2.5} Let $p>3$ be a prime. Then
$$S_p\Ls 12\e\cases-\f {p^2}{4x^2}\mod {p^3}\qq\q\qq\qq\qq\qq\t{if $p=x^2+4y^2\e
1\mod 4$,}
\\ \f{(p+1)^2}{2^{p-1}}
\b{(p-1)/2}{(p-3)/4}^2+\f {p^2}2\b{(p-1)/2}{(p-3)/4}^2E_{p-3}\mod
{p^3}\q\t{if $4\mid p-3$}.
\endcases$$
\endpro
Proof. Set $a=\f 12$. Then $\ap=\f{p+1}2$ and $a'=-\f 12$. For
$p=x^2+4y^2\e 1\mod 4$, since $2\nmid \ap$ we have
$$S_p\Ls 12\e \f{-\f 12\cdot \f 12\cdot p^2}{\f
14\b{(p-1)/2}{(p-1)/4}^2}\e -\f{p^2}{4x^2}\mod {p^3}$$ by Theorem
2.1. Now assume that $p\e 3\mod 4$. Then $\ap=2\cdot\f{p+1}4$. From
Theorem 2.1 and Lemma 2.4 we deduce that
$$\align S_p\Ls 12&\e\b{(p-1)/2}{(p+1)/4}^2\Big(1+p(H_{\f{p+1}2}+\qp
2) +\f{p^2}2\Big(-\qp 2^2+(H_{\f{p+1}2}+\qp 2)^2\\&\q-\f
14H_{\f{p+1}4}^{(2)}+\f 32H_{\f{p+1}2}^{(2)}\Big)\Big) \\&\e
\b{(p-1)/2}{(p+1)/4}^2\Big(1+p\Big(\f 2{p+1}-2\qp 2+p\qp 2^2+\qp
2\Big)
\\&\q+\f{p^2}2\Big(-\qp 2^2+\Big(\f 2{p+1}-2\qp 2+p\qp 2^2+\qp
2\Big)^2
\\&\q-\f 14\Big(\f{16}{(p+1)^2}-4E_{p-3}\Big)+\f
32\cdot\f 4{(p+1)^2}\Big)\Big)
\\&\e \b{(p-1)/2}{(p-3)/4}^2\Big(1+p(2-\qp 2)+p^2((1-\qp 2)^2
+\f 12E_{p-3})\Big)
\\&\e \f{(p+1)^2}{2^{p-1}}\b{(p-1)/2}{(p-3)/4}^2+\f
{p^2}2\b{(p-1)/2}{(p-3)/4}^2E_{p-3}
 \mod{p^3}.
\endalign$$
This proves the lemma.

 \pro{Theorem 2.2 ([S14, Conjecture 5.4])} Let
$p$ be a prime with $p>3$. Then
$$\align &\sum_{k=0}^{p-1}\f{\b{2k}k^3}{64^k(2k-1)}\\&\e
\cases -2x^2+p+\f{p^2}{4x^2} \mod {p^3}\q\qq\qq\qq\qq\t{if
$p=x^2+4y^2\e 1\mod 4$,}\\-\f
14\Big(2+p^2E_{p-3}\Big)\f{(p+1)^2}{2^{p-1}}
\b{(p-1)/2}{(p-3)/4}^2+\f {p^2}{2\b{(p-1)/2}{(p-3)/4}^2}\mod
{p^3}\q\t{if $4\mid p-3$.}
\endcases\endalign$$
\endpro
Proof. Note that
$$\b{\f 12}k\b{-\f 32}k=\f{1+2k}{1-2k}\b{-\f 12}k^2=\Big(-1-\f
2{2k-1}\Big)\b{2k}k^2\f 1{16^k}.\tag 2.6$$ We then have
$$2\sum_{k=0}^{p-1}\f{\b{2k}k^3}{64^k(2k-1)}+\sum_{k=0}^{p-1}\f{\b{2k}k^3}{64^k}
= -\sum_{k=0}^{p-1}\b{\f 12}k\b{-\f 32}k\b{2k}k\f 1{4^k}=-S_p\Ls
12.\tag 2.7$$ Now, applying (1.1) and Lemma 2.5 yields the result.

\section*{3. The congruence for
$\sum_{k=0}^{p-2}\b{a}k\b{-1-a}k{\b{2k}{k}}\f 1{4^k(k+1)}$ modulo
$p^3$}
\par Suppose that $F(a,k)$ and $R(a,k)$ satisfy $R(a,0)=0$ and
$$F(a,k)=F(a,k+1)R(a,k+1)-F(a,k)R(a,k)\q(k=0,1,2,\ldots).$$
For any positive integer $n$ we then have
$$\aligned \sum_{k=0}^{n-1}F(a,k)
&=\sum_{k=0}^{n-1}(F(a,k+1)R(a,k+1)-F(a,k)R(a,k))
\\&=F(a,n)R(a,n)-F(a,0)R(a,0)=F(a,n)R(a,n).\endaligned\tag 3.1$$
Using the WZ method and Zeilberger's algorithm one can easily find
$R(a,k)$ by Maple. See [PWZ].

\pro{Lemma 3.1} For any positive integer $n$ and real number
$a\not=0$ we have
$$\align &\sum_{k=0}^{n-1}\b{a}k\b{-1-a}k\b{2k}k\f
1{4^k(k+1)}\\&=S_n(a)+\f{a+1}aS_n(a+1)+F(a,n)R(a,n),\endalign$$
where
$$\align
&F(a,n)=\Big(\b{a}n\b{-1-a}n\big(\f
1{n+1}-1\big)-\f{a+1}a\b{a+1}n\b{-2-a}n\Big)\b{2n}n\f 1{4^n},
\\&R(a,n)=-\f{2n^2(n+1)}{n^2+2(a+1)^2n+(a+1)^2}.
\endalign$$
\endpro
Proof. It is easy to check that
$$R(a,0)=0,\q F(a,k)=F(a,k+1)R(a,k+1)-F(a,k)R(a,k).$$
Thus the result follows from (3.1).

\pro{Lemma 3.2} Let $p$ be an odd prime, $a\in\Bbb Z_p$, $a\not\e
0,-1\mod p$ and $a'=(a-\ap)/p$. Then
$$\f{\b a{p-1}\b{-1-a}{p-1}\b{2(p-1)}{p-1}}{4^{p-1}\cdot p}
\e -\f{a'(a'+1)}{a(a+1)}p^2\mod {p^3}.$$
\endpro
Proof. Since $\b a{p-1}=\f a{a+1-p}\b{a-1}{p-1}$ and
$\b{2(p-1)}{p-1}=\f p{2p-1}\b{2p-1}{p-1}$, we see that
$$\f{\b a{p-1}\b{-1-a}{p-1}\b{2(p-1)}{p-1}}{4^{p-1}\cdot p}
=\f a{(a+1-p)(2p-1)\cdot
4^{p-1}}\b{a-1}{p-1}\b{-1-a}{p-1}\b{2p-1}{p-1}.$$ By [S9, Lemma
2.2], $\b{a-1}{p-1}\b{-1-a}{p-1}\e \f{a'(a'+1)}{a^2}p^2\mod {p^3}$.
Thus, the result follows.

 \pro{Theorem 3.1} Let $p$ be an
odd prime, $a\in\Bbb Z_p$, $a\not\e 0,-1\mod p$ and $a'=(a-\ap)/p$.
Then
$$\align &\sum_{k=0}^{p-1}\b ak\b{-1-a}k\b{2k}k\f 1{4^k(k+1)}
+\f{a'(a'+1)}{a(a+1)}p^2
\\&\e \sum_{k=0}^{p-2}\b ak\b{-1-a}k\b{2k}k\f 1{4^k(k+1)}
\e S_p(a)+\f{a+1}aS_p(a+1)\mod {p^3}.\endalign$$ Thus, for $2\mid
\ap$  we have
$$\align &\sum_{k=0}^{p-2}\b ak\b{-1-a}k\b{2k}k\f 1{4^k(k+1)}
\\&\e
\b{(p-1)/2}{\ap/2}^2\Big(1+p\big((2a'+2)H_{\ap}-(2a'+1)H_{\f{\ap}2}
-2a'\qp2\big)
\\&\q+\f{p^2}2 \Big(2a'\qp 2^2+\big((2a'+2)H_{\ap}-(2a'+1)H_{\f{\ap}2}
-2a'\qp2\big)^2
\\&\q+\f 12(2{a'}^2-1)H_{\f{\ap}2}^{(2)}+2(1-{a'}^2)H_{\ap}^{(2)}\Big)\Big)
+p^2\f{a'(a'+1)}{a(a+1)}\b{(p-1)/2}{\ap/2}^{-2} \mod{p^3};
\endalign$$
for $2\nmid \ap$ and $\ap\not=p-2$ we have
$$\align &\sum_{k=0}^{p-2}\b ak\b{-1-a}k\b{2k}k\f 1{4^k(k+1)}
\\&\e\f {a+1}a
\b{(p-1)/2}{(\ap+1)/2}^2\Big(1+p\big((2a'+2)H_{\ap+1}-(2a'+1)H_{\f{\ap+1}2}
-2a'\qp2\big)
\\&\q+\f{p^2}2 \Big(2a'\qp 2^2+\big((2a'+2)H_{\ap+1}-(2a'+1)H_{\f{\ap+1}2}
-2a'\qp2\big)^2
\\&\q+\f 12(2{a'}^2-1)H_{\f{\ap+1}2}^{(2)}+2(1-{a'}^2)H_{\ap+1}^{(2)}\Big)
\Big) \\&\q+p^2\f{a'(a'+1)}{(a+1)^2}\cdot\b{(p-1)/2}{(\ap+1)/2}^{-2}
\mod{p^3}.
\endalign$$
\endpro
Proof. Since
$$\align&\b{a+1}p=\b{\ap+1+a'p}p\e a'\mod p,\\&
\b{-2-a}p=\b{p-2-\ap-(a'+1)p}p\e -(a'+1)\mod p,\endalign$$ taking
$n=p$ in Lemma 3.1 gives $$\align&  \sum_{k=0}^{p-1}\b
ak\b{-1-a}k\b{2k}k\f 1{4^k(k+1)} - S_p(a)-\f{a+1}aS_p(a+1)
\\&\e \f{a+1}a\b{a+1}p\b{-2-a}p\b{2p}p\f
1{4^p}\cdot\f{2p^2}{(a+1)^2}
\\&\e -\f{a'(a'+1)}{a(a+1)}p^2\mod {p^3}.
\endalign$$
Now, applying Lemma 3.2 and Theorem 2.1 yields the remaining part.

\pro{Theorem 3.2} Let $p$ be an odd prime. Then
$$\sum_{k=0}^{\f{p-1}2}\f{\b{2k}k^3}{64^k(k+1)}\e
\cases 4x^2-2p\mod {p^3}\qq\qq\t{if $p=x^2+4y^2\e 1\mod 4$,}
\\-\f{(p+1)^2}{2^{p-1}}
\b{(p-1)/2}{(p-3)/4}^2-p^2\Big(\b{(p-1)/2}{(p-3)/4}^{-2} \\\qq+\f
12\b{(p-1)/2}{(p-3)/4}^2E_{p-3}\Big)\mod {p^3}\q\qq\t{if $4\mid
p-3$}.
\endcases$$
\endpro
Proof. Note that $p\mid \b{2k}k$ for $\f p2<k<p$. Taking $a=-\f 12$
in Theorem 3.1 gives
$$\sum_{k=0}^{(p-1)/2}\f{\b{2k}k^3}{64^k(k+1)}
\e \sum_{k=0}^{p-2}\f{\b{2k}k^3}{64^k(k+1)}\e
\sum_{k=0}^{p-1}\f{\b{2k}k^3}{64^k}-S_p\Ls 12\mod {p^3}.$$ Now,
applying (1.1) and Lemma 2.5 yields the result.

\pro{Lemma 3.3} Let $p>3$ be a prime. Then
$$\align&H_{[\f p3]}\e -\f 32\qp 3+\f 34p\qp 3^2-p\Ls p3U_{p-3}\mod
{p^2},
\\&H_{[\f{2p}3]}\e -\f 32\qp 3+\f 34p\qp 3^2+2p\Ls p3U_{p-3}\mod
{p^2},
\\&H_{[\f p3]}^{(2)}\e -H_{[\f{2p}3]}^{(2)}\e 3\Ls p3U_{p-3}\mod p.
\endalign$$
\endpro
Proof. The first congruence was given in [S3, Theorem 3.2]. By [S3,
Theorem 3.2], $\sum_{k=1}^{[2p/3]}\f{(-1)^{k-1}}k\e 3p\sls
p3U_{p-3}\mod{p^2}$. Thus,
$$H_{[\f{2p}3]}=H_{[\f p3]}+\sum_{k=1}^{[2p/3]}\f{(-1)^{k-1}}k
\e -\f 32\qp 3+\f 34p\qp 3^2+2p\Ls p3U_{p-3}\mod{p^2}.$$ By [S3,
Theorem 3.3], $H_{[\f p3]}^{(2)}\e 3\ls p3U_{p-3}\mod p$. To
complete the proof, we note that
$$H_{[\f{2p}3]}^{(2)}=\sum_{k=1}^{p-1}\f 1{k^2}- \sum_{k=1}^{[p/3]}\f
1{(p-k)^2}\e -H_{[\f p3]}^{(2)}\mod p.$$

\pro{Theorem 3.3} Let $p$ be a prime with $p>3$. For $p\e 1\mod 3$
and so $p=x^2+3y^2$ we have
$$\sum_{k=0}^{p-2}\f{\b{2k}k^2\b{3k}k}{108^k(k+1)}
\e 4x^2-2p \mod {p^3}.$$ For $p\e 2\mod 3$ we have
$$\align
&\sum_{k=0}^{p-2}\f{\b{2k}k^2\b{3k}k}{108^k(k+1)}
\\&\e -2\b{(p-1)/2}{(p-5)/6}^2\Big(1+p\Big(2+\f 43\qp 2-\f 32\qp
3\Big) +p^2\Big(1+\f 83\qp 2+\f 29\qp 2^2\\&\q-3\qp 3-2\qp 2\qp
3+\f{15}8\qp 3^2+\f 34U_{p-3}\Big)\Big)-\f
12p^2\b{(p-1)/2}{(p-5)/6}^{-2} \mod {p^3}.
\endalign$$
\endpro
Proof. For $p=x^2+3y^2\e 1\mod 3$ we have $\langle \f
23\rangle_p=\f{p+2}3\e 1\mod 2$ and $\f 23-\langle \f
23\rangle_p=-\f p3$. Thus, from Theorem 2.1 and the well known fact
$\b{(p-1)/2}{(p-1)/6}\e 2x\sls x3\mod p$ (see [BEW]) we obtain
$$S_p\Ls 23\e\f{-\f 29p^2}{\f 49\b{(p-1)/2}{(p-1)/6}^2}
\e -\f{p^2}{8x^2}\mod {p^3}.$$ Taking $a=-\f 13$ in Theorem 3.1 and
then applying the above and (1.2) yields
$$\align\sum_{k=0}^{p-2}\f{\b{2k}k^2\b{3k}k}{108^k(k+1)}
&\e S_p\Big(-\f 13\Big)-2S_p\Ls 23 \e
4x^2-2p-\f{p^2}{4x^2}-2\Big(-\f{p^2}{8x^2}\Big)\\&=4x^2-2p
\mod{p^3}.\endalign$$

\par Now assume that $p\e 2\mod 3$. By Lemma 3.3,
$$\align\f 23H_{\f{2(p+1)}3}+\f 13H_{\f{p+1}3}
&\e\f 23\Big(\f 1{2(p+1)/3}-\f 32\qp 3+\f 34p\qp 3^2+2p\Ls
p3U_{p-3}\Big)
\\&\q+\f 13\Big(\f 1{(p+1)/3}-\f 32\qp 3+\f 34p\qp 3^2-p\Ls
p3U_{p-3}\Big)
\\&\e 2-\f 32\qp 3+p\Big(-2+\f 34\qp 3^2-U_{p-3}\Big)\mod {p^2}.\endalign$$
From Lemma 3.3 we also have
$$\align &H_{\f{p+1}3}^{(2)}=\f 9{(p+1)^2}+H_{[\f p3]}^{(2)}\e
9-3U_{p-3}\mod {p},
\\&H_{\f{2(p+1)}3}^{(2)}=\f 9{4(p+1)^2}+H_{[\f{2p}3]}^{(2)}\e \f 94
+3U_{p-3}\mod p.\endalign$$ Set $a=-\f 13$. Then $\ap=\f{2p-1}3$ and
$a'=-\f 23$. From Theorem 3.1 and the above,
$$\align &\sum_{k=0}^{p-2}\f{\b{2k}k^2\b{3k}k}{108^k(k+1)}
\\&\e -2\b{(p-1)/2}{(p+1)/3}^2
\Big(1+p\big(\f 23H_{\f{2(p+1)}3}+\f 13H_{\f{p+1}3}+\f 43\qp
2\big)\\&\q+\f{p^2}2\Big(-\f 43\qp 2^2+\big(\f 23H_{\f{2(p+1)}3}+\f
13H_{\f{p+1}3}+\f 43\qp 2\big)^2-\f 1{18}H_{\f{p+1}3}^{(2)}
+\f{10}9H_{\f{2(p+1)}3}^{(2)}\Big)\Big)\\&\q-\f
12p^2\b{(p-1)/2}{(p+1)/3}^{-2}
\\&\e -2\b{(p-1)/2}{(p-5)/6}^2\Big(1+p\Big(2+\f 43\qp 2-\f 32\qp
3\Big)+p^2\Big(-2+\f 34\qp 3^2-U_{p-3}\Big)
\\&\q+\f{p^2}2\Big(-\f 43\qp 2^2+\big(2+\f 43\qp 2-\f 32\qp
3\big)^2-\f 1{18}(9-3U_{p-3})+\f{10}9\big(\f
94+3U_{p-3}\big)\Big)\Big)\\&\q-\f 12p^2\b{(p-1)/2}{(p-5)/6}^{-2}
\\&\e -2\b{(p-1)/2}{(p-5)/6}^2\Big(1+p\Big(2+\f 43\qp 2-\f 32\qp
3\Big) +p^2\Big(1+\f 83\qp 2+\f 29\qp 2^2\\&\q-3\qp 3-2\qp 2\qp
3+\f{15}8\qp 3^2+\f 34U_{p-3}\Big)\Big)-\f
12p^2\b{(p-1)/2}{(p-5)/6}^{-2} \mod {p^3}.
\endalign$$
This completes the proof.

\pro{Theorem 3.4} Let $p>3$ be a prime. Then
$$\sum_{k=0}^{p-2}\f{\b{2k}k^2\b{4k}{2k}}{256^k(k+1)}
\e\cases 4x^2-2p\mod {p^3}&\t{if $p=x^2+2y^2\e 1,3\mod 8$,}
\\-\f 13R_2(p)\mod {p^2}&\t{if $p\e 5,7\mod 8$.}\endcases$$
\endpro
Proof. Taking $a=-\f 14$ or $-\f 34$ in Theorem 3.1 gives
$$\sum_{k=0}^{p-2}\f{\b{2k}k^2\b{4k}{2k}}{256^k(k+1)}\e
S_p\Big(-\f 14\Big)-3S_p\Ls 34 \e S_p\Big(-\f 34\Big)-\f 13S_p\Ls
14\mod {p^3}.$$ For $p=x^2+2y^2\e 1\mod 8$ we see that $\langle \f
34\rangle_p=\f{p+3}4\e 1\mod 2$ and $\f 34-\langle \f
34\rangle_p=-\f p4$. By Theorem 2.1 and the fact that
$\b{(p-1)/2}{(p-1)/8}\e 2(-1)^{\f{p-1}8+\f{x-1}2}x\mod p$ (see
[BEW]),
$$S_p\Ls 34\e \f{-\f 14(1-\f 14)p^2}{\sls
34^2\b{(p-1)/2}{(p-1)/8}^2}\e -\f{p^2}{12x^2}\mod {p^3}.$$ Recall
that $S_p(-\f 14)=\sum_{k=0}^{p-1}\f{\b{2k}k^2\b{4k}{2k}}{256^k}\e
4x^2-2p-\f{p^2}{4x^2}\mod {p^3}$. We then get
$$\align\sum_{k=0}^{p-2}\f{\b{2k}k^2\b{4k}{2k}}{256^k(k+1)}&\e
S_p\Big(-\f 14\Big)-3S_p\Ls 34 \e
4x^2-2p-\f{p^2}{4x^2}-3\Big(-\f{p^2}{12x^2}\Big) \\&=4x^2-2p\mod
{p^3}.
\endalign$$
\par For $p=x^2+2y^2\e 3\mod 8$ we see that $\langle \f
14\rangle_p=\f{p+1}4\e 1\mod 2$ and $\f 14-\langle \f
14\rangle_p=-\f p4$. By Theorem 2.1 and the fact that
$\b{(p-1)/2}{(p-3)/8}\e 2(-1)^{\f{p+5}8+\f{x-1}2}x\mod p$ (see
[BEW]), we deduce that
$$S_p\Ls 14\e \f{-\f 14(1-\f 14)p^2}{(-\f 14)^2\b{(p-1)/2}{(p-3)/8}^2}\e
-\f{3p^2}{4x^2}\mod {p^3}.$$ Recall that $S_p(-\f
34)=\sum_{k=0}^{p-1}\f{\b{2k}k^2\b{4k}{2k}}{256^k}\e
4x^2-2p-\f{p^2}{4x^2}\mod {p^3}$. We then get
$$\align\sum_{k=0}^{p-2}\f{\b{2k}k^2\b{4k}{2k}}{256^k(k+1)}&\e
S_p\Big(-\f 34\Big)-\f 13S_p\Ls 14 \e 4x^2-2p-\f{p^2}{4x^2}-\f
13\Big(-\f{3p^2}{4x^2}\Big) \\&=4x^2-2p\mod {p^3}.\endalign$$
\par For $p\e 5\mod 8$, taking $a=-\f 14$, $\ap=\f{p-1}4$ and
$a'=-\f 14$ in Theorem 3.1 and then applying the fact that $H_{[\f
p4]}\e -3\qp 2\mod p$ yields
$$\align&\sum_{k=0}^{p-2}\f{\b{2k}k^2\b{4k}{2k}}{256^k(k+1)}
\\&\e -3\b{(p-1)/2}{(p+3)/8}^2\Big(1+p\big(\f 32H_{\f{p+3}4}-\f
12H_{\f{p+3}8}+\f 12\qp 2\big)\Big)
\\&\e -3\Ls{3p+1}{p+3}^2\b{(p-1)/2}{(p-5)/8}^2\Big(1+p\big(\f
32\big(\f 4{p+3}-3\qp 2\big)-\f 12\big(\f 8{p+3}+H_{[\f p8]}\big)+\f
12\qp 2\big)\Big)
\\&\e \Big(-\f 13-\f{16}9p\Big)\b{(p-1)/2}{(p-5)/8}^2\Big(1+p\Big(\f
23-4\qp 2-\f 12H_{[\f p8]}\Big)\Big)\e -\f 13R_2(p)\mod {p^2}.
\endalign$$
\par For $p\e 7\mod 8$, taking $a=-\f 34$, $\ap=\f{p-3}4$ and
$a'=-\f 14$ in Theorem 3.1 and then applying the fact that $H_{[\f
p4]}\e -3\qp 2\mod p$ (see [S12,(2.4)]) yields
$$\align&\sum_{k=0}^{p-2}\f{\b{2k}k^2\b{4k}{2k}}{256^k(k+1)}
\\&\e -\f 13\b{(p-1)/2}{(p+1)/8}^2\Big(1+p\big(\f 32H_{\f{p+1}4}-\f
12H_{\f{p+1}8}+\f 12\qp 2\big)\Big)
\\&\e -\f 13\cdot 9\b{(p-1)/2}{(p-7)/8}^2\Big(1+p\big(\f 32\big(\f
4{p+1}-3\qp 2\big)-\f 12\big(\f 8{p+1}+H_{[\f p8]}\big)+\f 12\qp
2\big)\Big)
\\&\e-\f 13\cdot 9\b{(p-1)/2}{(p-7)/8}^2\Big(1+p\big(2-4\qp 2-\f
12H_{[\f p8]}\big)\Big)=-\f 13R_2(p)\mod {p^2}.
\endalign$$
\par Putting all the above together proves the theorem.

\pro{Theorem 3.5} Let $p$ be a prime with $p>3$. Then
$$\align &\sum_{k=0}^{p-2}\f{\b{2k}k\b{3k}k\b{6k}{3k}}{1728^k(k+1)}
\\&\e\cases
\sls p3(4x^2-2p)\mod {p^2}\qq\qq\qq\qq\q\ \t{if $p=x^2+4y^2\e 1\mod 4$,} \\
-\f 15\b{\f{p-1}2}{[\f p{12}]}^2\big(1+p\big(\f 25-3\qp 2-\f 52\qp
3-\f 23H_{[\f p{12}]}\big)\big)\mod {p^2}\q\t{if $12\mid p-7$,}
\\-5\b{\f{p-1}2}{[\f p{12}]}^2\big(1+p\big(2-3\qp 2-\f
52\qp 3-\f 23H_{[\f p{12}]}\big)\big)\mod {p^2}\q\t{if $12\mid
p-11$.}\endcases\endalign$$
\endpro
Proof. For $p=x^2+4y^2\e 1\mod 4$, taking $a=-\f 16$ in Theorem 3.1
and then applying Theorem 2.1 yields
$$\sum_{k=0}^{p-2}\f{\b{2k}k\b{3k}k\b{6k}{3k}}{1728^k(k+1)}
\e S_p\Big(-\f 16\Big)+\f {1-\f 16}{-\f 16}S_p\Big(-\f 16+1\Big)\e
S_p\Big(-\f 16\Big)\e \Ls p3(4x^2-2p)\mod {p^2}.$$
 For $p\e 7\mod
{12}$, taking $a=-\f 16$, $\ap=\f{p-1}6$ and $a'=-\f 16$ in Theorem
3.1 and then applying the fact that $H_{[\f p6]}\e -2\qp 2-\f 32\qp
3\mod p$ (see [S12,(2.4)]) yields
$$\align &\sum_{k=0}^{p-2}\f{\b{2k}k\b{3k}k\b{6k}{3k}}{1728^k(k+1)}
\\&\e -5\b{(p-1)/2}{(p+5)/12}^2\Big(1+p\Big(\f 53H_{\f{p+5}6}-\f
23H_{\f{p+5}{12}}+\f 13\qp 2\Big)\Big)
\\&\e -\f 15\b{(p-1)/2}{[p/12]}^2\Big(1+p\Big(\f 53\Big(\f
6{p+5}-2\qp 2-\f 32\qp 3\Big)-\f 23\Big(\f{12}{p+5}+H_{[\f
p{12}]}\Big)+\f 13\qp 2\Big)\Big)
\\&\e -\f 15\b{\f{p-1}2}{[\f p{12}]}^2\Big(1+p\big(\f 25-3\qp 2-\f
52\qp 3-\f 23H_{[\f p{12}]}\big)\Big)\mod {p^2}.
\endalign$$
For $p\e 11\mod {12}$, taking $a=-\f 56$, $\ap=\f{p-5}6$ and $a'=-\f
16$ in Theorem 3.1 and then applying the fact that $H_{[\f p6]}\e
-2\qp 2-\f 32\qp 3\mod p$ yields
$$\align &\sum_{k=0}^{p-2}\f{\b{2k}k\b{3k}k\b{6k}{3k}}{1728^k(k+1)}
\\&\e \f{1-\f 56}{-\f 56}\b{(p-1)/2}{(p+1)/12}^2\Big(1+p\Big(\f 53H_{\f{p+1}6}-\f
23H_{\f{p+1}{12}}+\f 13\qp 2\Big)\Big)
\\&\e -5\b{(p-1)/2}{[p/12]}^2\Big(1+p\Big(\f 53\Big(\f
6{p+1}-2\qp 2-\f 32\qp 3\Big)-\f 23\Big(\f{12}{p+1}+H_{[\f
p{12}]}\Big)+\f 13\qp 2\Big)\Big)
\\&\e -5\b{\f{p-1}2}{[\f p{12}]}^2\Big(1+p\big(2-3\qp 2-\f
52\qp 3-\f 23H_{[\f p{12}]}\big)\Big)\mod {p^2}.
\endalign$$
This completes the proof.
\par\q
\par{\bf Remark 3.1} Let $p>5$ be a prime. In the case $p\e 1\mod 4$, Theorem 3.5
is equivalent to a result due to Z.W. Sun [Su]. In [S14], the author
made a conjecture equivalent to
$$\Ls p3\sum_{k=0}^{p-2}\f{\b{2k}k\b{3k}k\b{6k}{3k}}{12^{3k}(k+1)}
\e\cases 4x^2-2p\mod {p^3}&\t{if $p=x^2+4y^2\e 1\mod 4$,}
\\\f 35R_1(p)\mod {p^2}&\t{if $p\e 3\mod 4$.}\endcases$$
For the conjectures concerning Theorems 3.2-3.4 see [S14,
Conjectures 5.4, 5.11 and 5.16].
\section*{4. The congruence for
$\sum_{k=0}^{p-1}\b{a}k\b{-1-a}k\b{2k}{k}\f{k}{4^k}$ modulo $p^3$}

\pro{Lemma 4.1} For any positive integer $n$ and real number $a$ we
have
$$\sum_{k=0}^{n-1}\b{a}k\b{-1-a}k\b{2k}k\f
{k}{4^k}-(a(a+1)S_n(a)-(a+1)^2S_n(a+1))=F(a,n)R(a,n),$$ where
$$\align
&F(a,n)=\Big(\b{a}n\b{-1-a}n(n-a(a+1))+(a+1)^2\b{a+1}n\b{-2-a}n\Big)\b{2n}n\f
1{4^n},
\\&R(a,n)=\f{2n^3}{n^2-2(a+1)^2n-(a+1)^2}.\endalign$$
\endpro
Proof. It is easy to check that $R(a,0)=0$ and
$F(a,k)=F(a,k+1)R(a,k+1)-F(a,k)R(a,k).$ Thus the result follows from
(3.1).
 \pro{Theorem 4.1} Let $p$ be an odd prime, $a\in\Bbb Z_p$,
$a\not\e -1\mod p$ and $a'=(a-\ap)/p$. Then
$$\align &\sum_{k=0}^{p-1}\b{a}k\b{-1-a}k\b{2k}k\f
{k}{4^k}\\&\e a(a+1)S_p(a)-(a+1)^2S_p(a+1)+\f{a'(a'+1)}{a+1}p^3\mod
{p^4}.\endalign$$ Moreover, for $a\not\e 0\mod p$,
$$\align&\sum_{k=0}^{p-1}\b{a}k\b{-1-a}k\b{2k}k\f
{k}{4^k}\\&\e
a(a+1)\Big(2S_p(a)-\sum_{k=0}^{p-2}\b{a}k\b{-1-a}k\b{2k}k\f
1{4^k(k+1)}\Big)\mod {p^3}.\endalign$$
\endpro
Proof. Since
$$\align &\b ap=\b{\ap+a'p}p\e a'\mod
p,\\&\b{-1-a}p=\b{p-1-\ap-(a'+1)p}p\e -(a'+1)\mod p,
\\&\b{a+1}p=\b{\ap+1+a'p}p\e a'\mod p,\\&
\b{-2-a}p=\b{p-2-\ap-(a'+1)p}p\e -(a'+1)\mod p,\endalign$$ putting
$n=p$ in Lemma 4.1 gives
$$\align &\sum_{k=0}^{p-1}\b{a}k\b{-1-a}k\b{2k}k\f
{k}{4^k}-(a(a+1)S_p(a)-(a+1)^2S_p(a+1))
\\&\e \Big(-a(a+1)\b
ap\b{-1-a}p+(a+1)^2\b{a+1}p\b{-2-a}p\Big)\b{2p}p\f 1{4^p}\cdot
\Big(-\f{2p^3}{(a+1)^2}\Big)
\\&\e (-a(a+1)a'(-a'-1)+(a+1)^2a'(-a'-1))\Big(-\f{p^3}{(a+1)^2}\Big)
=\f{a'(a'+1)}{a+1}p^3\mod {p^4}.\endalign$$ This together with
Theorem 3.1 yields
$$\align &\sum_{k=0}^{p-1}\b{a}k\b{-1-a}k\b{2k}k\f
{k}{4^k}+a(a+1)\sum_{k=0}^{p-2}\b{a}k\b{-1-a}k\b{2k}k\f 1{4^k(k+1)}
\\&\e 2a(a+1)S_p(a)\mod {p^3}.
\endalign$$
This completes the proof.

\pro{Theorem 4.2} Let $p$ be a prime with $p>3$. Then
$$\sum_{k=0}^{p-1}\f{k\b{2k}k^3}{64^k}\e
\cases -x^2+\f p2+\f{p^2}{8x^2} \mod {p^3}\q\qq\t{if $p=x^2+4y^2\e
1\mod 4$,}\\-\f{(p+1)^2}{2^{p+1}}
\b{(p-1)/2}{(p-3)/4}^2+\f{p^2}4\b{(p-1)/2}{(p-3)/4}^{-2}\\\qq-\f
{p^2}8\b{(p-1)/2}{(p-3)/4}^2E_{p-3}\mod {p^3}\q\qq\qq\t{if $4\mid
p-3$.}
\endcases$$
\endpro
Proof. Taking $a=-\f 12$ in Theorem 4.1 yields
$$ \sum_{k=0}^{p-1}\f{k\b{2k}k^3}{64^k}
\e -\f 14\Big(S_p\big(-\f 12\big)+S_p\ls 12\Big)\mod {p^3}.$$ Now
applying (1.1) and Lemma 2.5 gives the result.

\pro{Theorem 4.3} Let $p>3$ be a prime. For $p\e 1\mod 3$ and so
$p=x^2+3y^2$ we have
$$\sum_{k=0}^{p-1}\f{k\b{2k}k^2\b{3k}{k}}{108^k}
\e -\f 89x^2+\f 49p+\f{p^2}{9x^2}\mod {p^3}.$$ For $p\e 2\mod 3$ we
have
$$\align &\sum_{k=0}^{p-1}\f{k\b{2k}k^2\b{3k}{k}}{108^k}
\\&\e \f {p^2}9\b{(p-1)/2}{(p-5)/6}^{-2}-\f 49\b{(p-1)/2}{(p-5)/6}^2
\Big(1+p\Big(2+\f 43\qp 2-\f 32\qp 3\Big)
\\&\q+p^2\Big(1+\f 83\qp 2+\f 29\qp 2^2-3\qp 3-2\qp 2\qp
3+\f{15}8\qp 3^2+\f 34U_{p-3}\Big)\Big)\mod {p^3}.\endalign$$
\endpro
Proof. We first assume that $p=x^2+3y^2\e 1\mod 3$. Taking $a=-\f
13$ and $a'=-\f 13$ in Theorem 4.1 and then applying (1.2) and
Theorem 3.3 yields
$$\align\sum_{k=0}^{p-1}\f{k\b{2k}k^2\b{3k}{k}}{108^k}
&\e -\f 13\cdot\f
23\Big(2\big(4x^2-2p-\f{p^2}{4x^2}\big)-\big(4x^2-2p\big)\Big)
\\&=-\f 89x^2+\f 49p+\f{p^2}{9x^2}\mod {p^3}.\endalign$$
\par Now we assume that $p\e 2\mod 3$. Taking $a=-\f 13$ and $a'=-\f
23$ in Theorem 4.1 and then applying (1.2) and Theorem 3.3 yields
$$\align&\sum_{k=0}^{p-1}\f{k\b{2k}k^2\b{3k}{k}}{108^k}
\\&\e -\f 13\cdot \f
23\Big(-p^2\b{(p-1)/2}{(p-5)/6}^{-2}+2\b{(p-1)/2}{(p-5)/6}^2
\\&\q\times\Big(1+p\Big(2+\f 43\qp 2-\f 32\qp
3\Big) +p^2\Big(1+\f 83\qp 2+\f 29\qp 2^2\\&\q-3\qp 3-2\qp 2\qp
3+\f{15}8\qp 3^2+\f 34U_{p-3}\Big)\Big)+\f
{p^2}2\b{(p-1)/2}{(p-5)/6}^{-2}\Big)
\\&=\f {p^2}9\b{(p-1)/2}{(p-5)/6}^{-2}-\f 49\b{(p-1)/2}{(p-5)/6}^2
\Big(1+p\Big(2+\f 43\qp 2-\f 32\qp
3\Big)
\\&\q+p^2\Big(1+\f 83\qp 2+\f 29\qp 2^2-3\qp 3-2\qp 2\qp
3+\f{15}8\qp 3^2+\f 34U_{p-3}\Big)\Big)\mod {p^3}.
\endalign$$
This completes the proof.

\pro{Theorem 4.4} Let $p$ be an odd prime. Then
$$\sum_{k=0}^{p-1}\f{k\b{2k}k^2\b{4k}{2k}}{256^k} \e\cases
-\f 34x^2+\f 38p+\f{3p^2}{32x^2}\mod {p^3}& \t{if $p=x^2+2y^2\e
1,3\mod 8$,}
\\-\f 1{16}R_2(p) \mod {p^2}&\t{if $p\e 5,7\mod 8$.}
\endcases$$
\endpro
Proof. Taking $a=-\f 14$ in Theorem 4.1 yields
$$\sum_{k=0}^{p-1}\f{k\b{2k}k^2\b{4k}{2k}}{256^k} \e
-\f 14\cdot \f
34\Big(2\sum_{k=0}^{p-1}\f{\b{2k}k^2\b{4k}{2k}}{256^k}
-\sum_{k=0}^{p-2}\f{\b{2k}k^2\b{4k}{2k}}{256^k(k+1)}\Big)\mod{p^3}.$$
For $p=x^2+3y^2\e 1\mod 3$, applying (1.3) and Theorem 3.4 gives
$$\align\sum_{k=0}^{p-1}\f{k\b{2k}k^2\b{4k}{2k}}{256^k}
&\e -\f
3{16}\Big(2\big(4x^2-2p-\f{p^2}{4x^2}\big)-(4x^2-2p)\big)\\&=-\f
34x^2+\f 38p+\f{3p^2}{32x^2}\mod {p^3}.\endalign$$ For $p\e 2\mod
3$, from the above and Theorem 3.4 we deduce that
$$\sum_{k=0}^{p-1}\f{k\b{2k}k^2\b{4k}{2k}}{256^k}\e -\f
3{16}\Big(2\cdot 0+\f 13R_2(p)\Big)=-\f 1{16}R_2(p) \mod{p^2}.$$ The
proof is now complete.

\pro{Theorem 4.5} Let $p>3$ be a prime. Then
$$\align &\sum_{k=0}^{p-1}\f{k\b{2k}k\b{3k}k\b{6k}{3k}}{12^{3k}}
\\&\e \cases\sls p3(-\f 59x^2+\f{5}{18}p)\mod {p^2}\qq\qq\qq\qq\t{if $p=x^2+4y^2\e
1\mod 4$,}
\\-\f 1{36}\b{\f{p-1}2}{[\f p{12}]}^2\big(1+p\big(\f 25-3\qp 2-\f
52\qp 3-\f 23H_{[\f p{12}]}\big)\big)\mod {p^2}\q\t{if $12\mid
p-7$,}
\\-\f{25}{36}\b{\f{p-1}2}{[\f p{12}]}^2\big(1+p\big(2-3\qp 2-\f
52\qp 3-\f 23H_{[\f p{12}]}\big)\big)\mod{p^2}\q\t{if $12\mid
p-11$.}
\endcases\endalign$$
\endpro
Proof. Taking $a=-\f 16$ in Theorem 4.1 we get
$$\sum_{k=0}^{p-1}\f{k\b{2k}k\b{3k}k\b{6k}{3k}}{12^{3k}}
\e -\f
5{36}\Big(2\sum_{k=0}^{p-1}\f{\b{2k}k\b{3k}k\b{6k}{3k}}{12^{3k}}
-\sum_{k=0}^{p-2}\f{\b{2k}k\b{3k}k\b{6k}{3k}}{12^{3k}(k+1)}\Big)\mod
{p^3}.\tag 4.1$$ For $p=x^2+4y^2\e 1\mod 4$, applying Theorem 3.5 we
get
$$\align&\sum_{k=0}^{p-1}\f{k\b{2k}k\b{3k}k\b{6k}{3k}}{12^{3k}}
\\&
\e -\f 5{36}\Big(2\ls p3(4x^2-2p)-\ls p3(4x^2-2p)\Big)=\Ls
p3\Big(-\f 59x^2+\f 5{18}p\Big)\mod {p^2}.
\endalign$$
For $p\e 3\mod 4$ we have
$\sum_{k=0}^{p-1}\f{\b{2k}k\b{3k}k\b{6k}{3k}}{12^{3k}}\e 0\mod
{p^2}$. Thus, from (4.1) and Theorem 3.5 we deduce  the result in
this case.
\par{\bf Remark 4.1} In [S4] the author established the congruence
for $\sum_{k=0}^{p-1}\f{k\b{2k}k^3}{64^k}$ modulo $p^2$, where $p$
is a prime with $p>5$. For the conjectures concerning Theorems
4.3-4.5 see [S13, Conjectures 2.7, 2.11 and 2.16].
\section*{5. The congruence for
$\sum_{k=0}^{p-1}\b{a}k\b{-1-a}k\b{2k}{k}\f 1{4^k(2k-1)}$ modulo
$p^3$}

\pro{Lemma 5.1} For any positive integer $n$ and real number $a$ we
have
$$\align &\sum_{k=0}^{n-1}\b{a}k\b{-1-a}k\b{2k}k\f
1{4^k(2k-1)}\\&=-(2a^2+2a+1)S_n(a)-2(a+1)^2S_n(a+1)+F(a,n)R(a,n),\endalign$$
where
$$\align
&F(a,n)=\Big(\b{a}n\b{-1-a}n\big(\f
1{2n-1}+(2a^2+2a+1)\big)\\&\qq\qq+2(a+1)^2\b{a+1}n\b{-2-a}n\Big)\b{2n}n\f
1{4^n},
\\&R(a,n)=-\f{2n^3}{n^2+2(a+1)^2n-(a+1)^2}.\endalign$$
\endpro
Proof. It is easy to check that $R(a,0)=0$ and
$F(a,k)=F(a,k+1)R(a,k+1)-F(a,k)R(a,k).$ Thus the result follows from
(3.1). \pro{Theorem 5.1} Let $p$ be an odd prime, $a\in\Bbb Z_p$,
$a\not\e -1\mod p$ and $a'=(a-\ap)/p$. Then
$$\align &\sum_{k=0}^{p-1}\b ak\b{-1-a}k\b{2k}k\f 1{4^k(2k-1)}
\\&\e -(2a^2+2a+1)S_p(a)-2(a+1)^2S_p(a+1)-2(2a+1)\f{a'(a'+1)}{a+1}p^3
\mod {p^4}.
\endalign$$
\endpro
Proof. By the proof of Theorem 4.1,
$$\b ap\b{-1-a}p\e \b{a+1}p\b{-2-a}p\e -a'(a'+1)\mod p.$$
Taking $n=p$ in Lemma 5.1 and then applying the above yields
$$\align &\sum_{k=0}^{p-1}\b ak\b{-1-a}k\b{2k}k\f 1{4^k(2k-1)}
\\&+(2a^2+2a+1)S_p(a)+2(a+1)^2S_p(a+1)
\\&\e -a'(a'+1)\Big(\f 1{2p-1}+2a^2+2a+1+2(a+1)^2\Big)
\b{2p}p\f 1{4^p}\cdot\f{2p^3}{(a+1)^2}
\\&\e -a'(a'+1)(2a(a+1)+2(a+1)^2)\b{2p-1}{p-1}\f
1{4^{p-1}}\cdot\f{p^3}{(a+1)^2}
\\&\e -2(2a+1)\f{a'(a'+1)}{a+1}p^3
\mod {p^4}.
\endalign$$
This proves the theorem.

\pro{Corollary 5.1} Let $p$ be an odd prime, $a\in\Bbb Z_p$ and
$a\not\e 0,-1\mod p$. Then
$$\align &\sum_{k=0}^{p-1}\b{a}k\b{-1-a}k\b{2k}k\f
1{4^k(2k-1)}\\&\e
-S_p(a)-2a(a+1)\sum_{k=0}^{p-2}\b{a}k\b{-1-a}k\b{2k}k\f
1{4^k(k+1)}\mod {p^3}.
\endalign$$
\endpro
 Proof. This is immediate from Theorems 3.1 and 5.1.
 \pro{Theorem
5.2} Let $p$ be an odd prime. Then
$$\align &\sum_{k=0}^{p-1}\f{\b{2k}k^3}{64^k(2k-1)^2}\\&\e \cases 2x^2-p-\f{p^2}{2x^2}
\mod {p^3}\qq\qq\qq\q\qq\t{if $p=x^2+4y^2\e 1\mod 4$,}
\\\f{3(p+1)^2}{2^p}\b{\f{p-1}2}{\f{p-3}4}^2+\f
34p^2\b{\f{p-1}2}{\f{p-3}4}^2E_{p-3}
-\f{p^2}2\b{\f{p-1}2}{\f{p-3}4}^{-2}\mod{p^3}\q\t{if $4\mid p-3$.}
\endcases\endalign$$
\endpro
Proof. Taking $a=\f 12$ in Theorem 5.1 and then applying Lemma 2.1
gives
$$\align&\sum_{k=0}^{p-1}\b{\f 12}k\b{-\f 32}k\b{2k}k\f 1{4^k(2k-1)}
\\&\e -\Big(2\cdot \f 14+2\cdot \f 12+1\Big)S_p\Ls 12-2\Big(\f
12+1\Big)^2S_p\Big(\f 12+1\Big) \\&\e -\f 52S_p\Ls 12-\f
12S_p\Big(-\f 12\Big)\mod {p^3}.\endalign$$ Since
$\b{1/2}k\b{-3/2}k=-(1+\f 2{2k-1})\b{2k}k^2\f 1{16^k}$ by (2.6),
applying (2.7) and the above we see that
$$\align 2\sum_{k=0}^{p-1}\f{\b{2k}k^3}{64^k(2k-1)^2}
&=-\sum_{k=0}^{p-1}\f{\b{2k}k^3}{64^k(2k-1)}-\sum_{k=0}^{p-1}\b{\f
12}k\b{-\f 32}k\b{2k}k\f 1{4^k(2k-1)}
\\&\e \f 12\Big(S_p\Ls 12+S_p\Big(-\f 12\Big)\Big)+
\f 52S_p\Ls 12+\f 12S_p\Big(-\f 12\Big)\\&=3S_p\Ls 12+S_p\Big(-\f
12\Big)\mod {p^3}.\endalign$$ Now, applying (1.1) and Lemma 2.5
yields the result.

 \pro{Theorem
5.3} Let $p>3$ be a prime. For $p\e 1\mod 3$ and so $p=x^2+3y^2$ we
have
$$\sum_{k=0}^{p-1}\f{\b{2k}k^2\b{3k}k}{108^k(2k-1)}\e
-\f{20}9x^2+\f{10}9p+\f {p^2}{4x^2}\mod {p^3}.$$ For $p\e 2\mod 3$
we have
$$\align&\sum_{k=0}^{p-1}\f{\b{2k}k^2\b{3k}k}{108^k(2k-1)}
\\&\e-\f 89\b{\f{p-1}2}{\f{p-5}6}^2\Big(1+p\Big(2+\f 43\qp 2-\f 32\qp
3\Big) +p^2\Big(1+\f 83\qp 2+\f 29\qp 2^2\\&\q-3\qp 3-2\qp 2\qp
3+\f{15}8\qp 3^2+\f 34U_{p-3}\Big)\Big)+\f
5{18}p^2\b{\f{p-1}2}{\f{p-5}6}^{-2} \mod {p^3}.
\endalign$$
\endpro
Proof. Taking $a=-\f 13$ in Corollary 5.1 and noting that $a'=-\f
13$ or $-\f 23$ gives
$$\sum_{k=0}^{p-1}\f{\b{2k}k^2\b{3k}k}{108^k(2k-1)}\e
-\sum_{k=0}^{p-1}\f{\b{2k}k^2\b{3k}k}{108^k}+\f
49\sum_{k=0}^{p-2}\f{\b{2k}k^2\b{3k}k}{108^k(k+1)}\mod {p^3}.$$ Now
applying (1.2) and Theorem 3.3 yields the result.

\pro{Theorem 5.4} Let $p$ be an odd prime. Then
$$\sum_{k=0}^{p-1}\f{\b{2k}k^2\b{4k}{2k}}{256^k(2k-1)}
\e\cases -\f 52x^2+\f 54p+\f{p^2}{4x^2}\mod {p^3}&\t{if
$p=x^2+2y^2\e 1,3\mod 8$,}
\\-\f 18R_2(p)\mod {p^2}&\t{if $p\e 5,7\mod 8$.}
\endcases$$
\endpro
Proof. Taking $a=-\f 14$ in Corollary 5.1 and noting that $a'=-\f
14$ or $-\f 34$ gives
$$\sum_{k=0}^{p-1}\f{\b{2k}k^2\b{4k}{2k}}{256^k(2k-1)}
\e -\sum_{k=0}^{p-1}\f{\b{2k}k^2\b{4k}{2k}}{256^k}+\f
38\sum_{k=0}^{p-2}\f{\b{2k}k^2\b{4k}{2k}}{256^k(k+1)}\mod {p^3}.$$
Now applying (1.3) and Theorem 3.4 yields the result.

\pro{Theorem 5.5} Let $p>3$ be a prime. Then
$$\align &\sum_{k=0}^{p-1}\f{\b{2k}k\b{3k}k\b{6k}{3k}}{12^{3k}(2k-1)}
\\&\e\cases \sls p3(-\f {26}9x^2+\f {13}9p)\mod {p^2}\qq\qq\qq\qq\t{if
$p=x^2+4y^2\e 1\mod 4$,}
\\-\f 1{18}\b{\f{p-1}2}{[\f p{12}]}^2\big(1+p\big(\f 25-3\qp 2-\f
52\qp 3-\f 23H_{[\f p{12}]}\big)\big)\mod {p^2}\q\t{if $12\mid
p-7$,}
\\-\f {25}{18}\b{\f{p-1}2}{[\f p{12}]}^2\big(1+p\big(2-3\qp 2-\f
52\qp 3-\f 23H_{[\f p{12}]}\big)\big)\mod {p^2}\q\t{if $12\mid
p-11$.}\endcases\endalign$$
\endpro
Proof. Taking $a=-\f 16$ in Corollary 5.1 and noting that $a'=-\f
16$ or $-\f 56$ gives
$$\sum_{k=0}^{p-1}\f{\b{2k}k\b{3k}k\b{6k}{3k}}{12^{3k}(2k-1)}
\e -\sum_{k=0}^{p-1}\f{\b{2k}k\b{3k}k\b{6k}{3k}}{12^{3k}} +\f 5{18}
\sum_{k=0}^{p-2}\f{\b{2k}k\b{3k}k\b{6k}{3k}}{12^{3k}(k+1)}\mod
{p^3}.$$ Now, applying Theorem 3.5 yields the result.
\par{\bf Remark 5.1} For the conjectures concerning Theorems 5.2-5.5 see [S14,
Conjectures 5.4, 5.7, 5.11 and 5.16].
\section*{6. The congruence for
$\sum_{k=0}^{p-1}\b{a}k\b{-1-a}k\b{2k}{k}\f{k^2}{4^k}$ modulo $p^3$}
\pro{Lemma 6.1} For any positive integer $n$ and real number $a$ we
have
$$\align &\sum_{k=0}^{n-1}\b{a}k\b{-1-a}k\b{2k}k\f
1{4^k}\big(3k^2-(2a^2+2a-1)k-a(a+1)\big)
 \\&=2n^3\b{a}n\b{-1-a}n\b{2n}n\f
1{4^n}.\endalign$$
\endpro
Proof. Set
$$\align &F(a,k)=\b{a}k\b{-1-a}k\b{2k}k\f
1{4^k}\big(3k^2-(2a^2+2a-1)k-a(a+1)\big),
\\&G(a,k)=2k^3\b{a}k\b{-1-a}k\b{2k}k\f
1{4^k}.\endalign$$  It is easy to check that $G(a,0)=0$ and
$F(a,k)=G(a,k+1)-G(a,k).$ Thus
$\sum_{k=0}^{n-1}F(a,k)=G(a,n)-G(a,0)=G(a,n)$.
 \pro{Theorem
6.1} Let $p$ be an odd prime, $a\in\Bbb Z_p$ and $a'=(a-\ap)/p$.
Then
$$\align &\sum_{k=0}^{p-1}\b{a}k\b{-1-a}k\b{2k}k\f
1{4^k}\big(3k^2-(2a^2+2a-1)k-a(a+1)\big)\\&\e -a'(a'+1)p^3\mod
{p^4}.\endalign$$
\endpro
Proof. Taking $n=p$ in Lemma 6.1 and noting that $\b ap\e a'\mod p$,
 $\b{-1-a}p\e -(a'+1)\mod p$ and $\b{2p}p\f 1{4^p}=\f 1{2\cdot 4^{p-1}}
 \b{2p-1}{p-1}\e \f 12\mod p$ yields the result.
\par\q
\par From Theorems 4.1 and 6.1 we deduce the following result.
\pro{Theorem 6.2} Let $p$ be an odd prime, $a\in\Bbb Z_p$, $a\not\e
0,-1\mod p$ and $a'=(a-\ap)/p$. Then
$$\align \sum_{k=0}^{p-1}\b{a}k\b{-1-a}k\b{2k}k\f
{k^2}{4^k} &\e \f
23a^2(a+1)^2S_p(a)-\f{2a(a+1)-1}3(a+1)^2S_p(a+1)\\&\q+
\f{(2a^2+a-2)a'(a'+1)}{3(a+1)}p^3\mod {p^4}.
\endalign$$
\endpro

\pro{Theorem 6.3} Let $p$ be a prime with $p>3$. Then
$$\sum_{k=0}^{p-1}\f{k^2\b{2k}k^3}{64^k}\e
\cases \f 16x^2-\f p{12}-\f{p^2}{24x^2} \mod {p^3}\q\qq\t{if
$p=x^2+4y^2\e 1\mod 4$,}\\\f{(p+1)^2}{2^{p+2}}
\b{(p-1)/2}{(p-3)/4}^2-\f{p^2}{24}\b{(p-1)/2}{(p-3)/4}^{-2}\\\qq+\f
{p^2}{16}\b{(p-1)/2}{(p-3)/4}^2E_{p-3}\mod {p^3}\q\qq\qq\t{if $4\mid
p-3$.}
\endcases$$
\endpro
Proof. Taking $a=-\f 12$ in Theorem 6.1 gives
$$\sum_{k=0}^{p-1}\f{k^2\b{2k}k^3}{64^k}\e
-\f 12\sum_{k=0}^{p-1}\f{k\b{2k}k^3}{64^k}-\f
1{12}\sum_{k=0}^{p-1}\f{\b{2k}k^3}{64^k}\mod {p^3}.$$ Now, appealing
to (1.1) and Theorem 4.2 we deduce the result.

\pro{Theorem 6.4} Let $p>3$ be a prime. For $p\e 1\mod 3$ and so
$p=x^2+3y^2$ we have
$$\sum_{k=0}^{p-1}\f{k^2\b{2k}k^2\b{3k}{k}}{108^k}
\e \f {32}{243}x^2-\f {16}{243}p-\f{17p^2}{486x^2}\mod {p^3}.$$ For
$p\e 2\mod 3$ we have
$$\align &\sum_{k=0}^{p-1}\f{k^2\b{2k}k^2\b{3k}{k}}{108^k}
\\&\e -\f {4}{243}p^2\b{(p-1)/2}{(p-5)/6}^{-2}+\f {52}{243}\b{(p-1)/2}{(p-5)/6}^2
\Big(1+p\Big(2+\f 43\qp 2-\f 32\qp 3\Big)
\\&\q+p^2\Big(1+\f 83\qp 2+\f 29\qp 2^2-3\qp 3-2\qp 2\qp
3+\f{15}8\qp 3^2+\f 34U_{p-3}\Big)\Big)\mod {p^3}.\endalign$$
\endpro
Proof. Taking $a=-\f 13$ in Theorem 6.1 yields
$$\sum_{k=0}^{p-1}\f{k^2\b{2k}k^2\b{3k}{k}}{108^k}\e -\f{13}{27}
\sum_{k=0}^{p-1}\f{k\b{2k}k^2\b{3k}{k}}{108^k}-\f
2{27}\sum_{k=0}^{p-1}\f{\b{2k}k^2\b{3k}{k}}{108^k}\mod {p^3}.$$ Now
applying (1.2) and Theorem 4.3 deduces the result.

\pro{Theorem 6.5} Let $p$ be a prime with $p>3$. Then
$$\sum_{k=0}^{p-1}\f{k^2\b{2k}k^2\b{4k}{2k}}{256^k} \e\cases \f
3{32}x^2-\f 3{64}p-\f{7p^2}{256x^2}\mod {p^3}& \t{if $p=x^2+2y^2\e
1,3\mod 8$,}
\\\f {11}{384}R_2(p) \mod {p^2}& \t{if $p\e 5,7\mod 8$.}
\endcases$$
\endpro
Proof. Taking $a=-\f 14$ in Theorem 6.1 yields
$$\sum_{k=0}^{p-1}\f{k^2\b{2k}k^2\b{4k}{2k}}{256^k}\e -\f{11}{24}
\sum_{k=0}^{p-1}\f{k\b{2k}k^2\b{4k}{2k}}{256^k}-\f
1{16}\sum_{k=0}^{p-1}\f{\b{2k}k^2\b{4k}{2k}}{256^k}\mod {p^3}.$$
This together with (1.3) and Theorem 4.4 deduces the result.

\pro{Theorem 6.6} Let $p>3$ be a prime. Then
$$\align &\sum_{k=0}^{p-1}\f{k^2\b{2k}k\b{3k}k\b{6k}{3k}}{12^{3k}}
\\&\e \cases\sls p3(\f {25}{486}x^2-\f{25}{972}p)\mod {p^2}\qq\qq\qq\qq\t{if $p=x^2+4y^2\e
1\mod 4$,}
\\\f {23}{1944}\b{\f{p-1}2}{[\f p{12}]}^2\big(1+p\big(\f 25-3\qp 2-\f
52\qp 3-\f 23H_{[\f p{12}]}\big)\big)\mod {p^2}\q\t{if $12\mid
p-7$,}
\\\f{575}{1944}\b{\f{p-1}2}{[\f p{12}]}^2\big(1+p\big(2-3\qp 2-\f
52\qp 3-\f 23H_{[\f p{12}]}\big)\big)\mod{p^2}\q\t{if $12\mid
p-11$.}
\endcases\endalign$$
\endpro
Proof. Taking $a=-\f 16$ in Theorem 6.1 we get
$$\sum_{k=0}^{p-1}\f{k^2\b{2k}k\b{3k}k\b{6k}{3k}}{12^{3k}}
\e -\f{23}{54}\sum_{k=0}^{p-1}\f{k\b{2k}k\b{3k}k\b{6k}{3k}}{12^{3k}}
-\f 5{108}\sum_{k=0}^{p-1}\f{\b{2k}k\b{3k}k\b{6k}{3k}}{12^{3k}}\mod
{p^3}.$$ This together with Theorem 4.5 yields the result.

\par{\bf Remark 6.1} In [S4] the author established the congruence
for $\sum_{k=0}^{p-1}\f{k^2\b{2k}k^3}{64^k}$ modulo $p^2$, where $p$
is a prime with $p>5$. For the conjectures concerning Theorems
6.4-6.6 see [S13, Conjectures 2.7, 2.11 and 2.16].

\section*{7. The congruence for
$\sum_{k=0}^{p-1}\b{a}k\b{-1-a}k{\b{2k}{k}}\f{k^3}{4^k}$ modulo
$p^3$} \pro{Lemma 7.1} For any positive integer $n$ and real number
$a$ we have
$$\align&\sum_{k=0}^{n-1}\b{a}k\b{-1-a}k\f{\b{2k}k}{4^k}
\Big(15k^3-(4a^2(a+1)^2-a(a+1)+1)k-a(a+1)(2a(a+1)-1)\Big)
\\&=2n^3(3n+2a(a+1)-4)\b{a}n\b{-1-a}n\b{2n}n\f
1{4^n}.\endalign$$
\endpro
Proof. Set
$$\align &F(a,k)=\b{a}k\b{-1-a}k\b{2k}k\f 1{4^k}
\big(15k^3-(4a^2(a+1)^2-a(a+1)+1)k\\&\qq\qq\q-a(a+1)(2a(a+1)-1)\big),
\\&G(a,k)=2k^3(3k+2a(a+1)-4)\b{a}k\b{-1-a}k\b{2k}k\f
1{4^k}.\endalign$$ It is easy to check that $G(a,0)=0$ and
$F(a,k)=G(a,k+1)-G(a,k).$ Thus
$\sum_{k=0}^{n-1}F(a,k)=G(a,n)-G(a,0)=G(a,n)$.

 \pro{Theorem 7.1} Let $p$ be an odd prime, $a\in\Bbb Z_p$
 and $a'=(a-\ap)/p$. Then
$$\align &\sum_{k=0}^{p-1}k^3\b{a}k\b{-1-a}k\f{\b{2k}k}{4^k}\\&\q\times\big(15k^3-
(4a^2(a+1)^2-a(a+1)+1)k-a(a+1)(2a(a+1)-1)\big)
\\&\e (4-2a(a+1))a'(a'+1)p^3\mod
{p^4}.\endalign$$
\endpro
Proof. Taking $n=p$ in Lemma 7.1 and noting that $\b ap\e a'\mod p$,
 $\b{-1-a}p\e -(a'+1)\mod p$ and $\b{2p}p\f 1{4^p}=\f 1{2\cdot 4^{p-1}}
 \b{2p-1}{p-1}\e \f 12\mod p$ yields the result.

 \pro{Theorem 7.2} Let $p>5$ be a prime, $a\in\Bbb Z_p$, $a\not\e
 -1\mod p$
 and $a'=(a-\ap)/p$. Then
$$\align &\sum_{k=0}^{p-1}k^3\b{a}k\b{-1-a}k\b{2k}k\f 1{4^k}
\\&\e \f {a^2(a+1)^2(2a+1)^2}{15}S_p(a)-\f{(a+1)^2(4a^2(a+1)^2-a(a+1)+1)}{15}
S_p(a+1)\\&\q+\f{4a^4+6a^3-a^2+a+5}{15(a+1)}a'(a'+1)p^3\mod {p^4}.
\endalign$$
\endpro
Proof. This is immediate from Theorems 4.1 and 7.1.

\pro{Theorem 7.3} Let $p$ be a prime with $p>5$. Then
$$\sum_{k=0}^{p-1}\f{k^3\b{2k}k^3}{64^k}\e
\cases \f{p^2}{160x^2} \mod {p^3}\q\qq\qq\qq\t{if $p=x^2+4y^2\e
1\mod 4$,}\\-\f{(p+1)^2}{40\cdot 2^{p-1}}
\b{\f{p-1}2}{\f{p-3}4}^2-\f
{p^2}{80}\b{\f{p-1}2}{\f{p-3}4}^2E_{p-3}\mod {p^3}\q\t{if $4\mid
p-3$.}
\endcases$$
\endpro
Proof. Putting $a=-\f 12$ in Theorem 7.1 gives
$$\sum_{k=0}^{p-1}\f{k^3\b{2k}k^3}{64^k}\e
\f 1{10}\sum_{k=0}^{p-1}\f{k\b{2k}k^3}{64^k}+\f
1{40}\sum_{k=0}^{p-1}\f{\b{2k}k^3}{64^k}\mod {p^3}.$$ Now applying
(1.1) and Theorem 4.2 yields the result.
\par {\bf Remark 7.1}  Theorem 7.3 can be easily deduced
from a result due to Tauraso [T3],
 although he did not give the details of its proof.

\pro{Theorem 7.4} Let $p>5$ be a prime. For $p\e 1\mod 3$ and so
$p=x^2+3y^2$ we have
$$\sum_{k=0}^{p-1}\f{k^3\b{2k}k^2\b{3k}{k}}{108^k}
\e  \f 1{10935}\Big(16x^2-8p+\f{113p^2}{2x^2}\Big)\mod {p^3}.$$ For
$p\e 2\mod 3$ we have
$$\align &\sum_{k=0}^{p-1}\f{k^3\b{2k}k^2\b{3k}{k}}{108^k}
\\&\e -\f {2}{10935}p^2\b{(p-1)/2}{(p-5)/6}^{-2}-\f {92}{2187}\b{(p-1)/2}{(p-5)/6}^2
\Big(1+p\Big(2+\f 43\qp 2-\f 32\qp 3\Big)
\\&\q+p^2\Big(1+\f 83\qp 2+\f 29\qp 2^2-3\qp 3-2\qp 2\qp
3+\f{15}8\qp 3^2+\f 34U_{p-3}\Big)\Big)\mod {p^3}.\endalign$$
\endpro
Proof. Putting $a=-\f 13$ in Theorem 7.1 we obtain
$$\sum_{k=0}^{p-1}\f{k^3\b{2k}k^2\b{3k}{k}}{108^k}
\e \f{23}{243}\sum_{k=0}^{p-1}\f{k\b{2k}k^2\b{3k}{k}}{108^k}
+\f{26}{1215}\sum_{k=0}^{p-1}\f{\b{2k}k^2\b{3k}{k}}{108^k}\mod
{p^3}.$$ This together with (1.2) and Theorem 4.3 yields the result.

\pro{Theorem 7.5} Let $p>3$ be a prime. Then
$$\sum_{k=0}^{p-1}\f{k^3\b{2k}k^2\b{4k}{2k}}{256^k} \e\cases \f
{3x^2}{1280}-\f {3p}{2560}+\f{41p^2}{10240x^2}\mod {p^3}\ \; \t{if
$p=x^2+2y^2\e 1,3\mod 8$,}
\\-\f {17}{3072}R_2(p) \mod {p^2}\qq\qq\q\ \; \t{if $p\e 5,7\mod 8$.}
\endcases$$
\endpro
Proof. Putting $a=-\f 14$ in Theorem 7.1 we see that
$$\sum_{k=0}^{p-1}\f{k^3\b{2k}k^2\b{4k}{2k}}{256^k}
\e \f{17}{192}\sum_{k=0}^{p-1}\f{k\b{2k}k^2\b{4k}{2k}}{256^k}
+\f{11}{640}\sum_{k=0}^{p-1}\f{\b{2k}k^2\b{4k}{2k}}{256^k}\mod{p^3}.$$
Now applying (1.3) and Theorem 4.4 yields the result.

\pro{Theorem 7.6} Let $p>5$ be a prime. Then
$$\align&\sum_{k=0}^{p-1}\f{k^3\b{2k}k\b{3k}{k}\b{6k}{3k}}{12^{3k}}
\\&\q\e\cases \sls p3(\f{5}{2187}x^2-\f{5}{4374}p)\mod
{p^2}\qq\qq\qq\qq\t{if $p=x^2+4y^2\e 1\mod 4$,}
\\-\f {17}{6912}\b{\f{p-1}2}{[\f p{12}]}^2\big(1+p\big(\f 25-3\qp 2-\f
52\qp 3-\f 23H_{[\f p{12}]}\big)\big)\mod {p^2}\q\t{if $12\mid
p-7$,}
\\-\f{425}{6912}\b{\f{p-1}2}{[\f p{12}]}^2\big(1+p\big(2-3\qp 2-\f
52\qp 3-\f 23H_{[\f p{12}]}\big)\big)\mod{p^2}\q\t{if $12\mid
p-11$.}
\endcases
\endalign$$
\endpro
Proof. Putting $a=-\f 16$ in Theorem 7.1 we deduce that
$$\sum_{k=0}^{p-1}\f{k^3\b{2k}k\b{3k}{k}\b{6k}{3k}}{12^{3k}}
\e
\f{197}{2430}\sum_{k=0}^{p-1}\f{k\b{2k}k\b{3k}{k}\b{6k}{3k}}{12^{3k}}
+\f{23}{1944}\sum_{k=0}^{p-1}\f{\b{2k}k\b{3k}{k}\b{6k}{3k}}{12^{3k}}\mod
{p^3}.$$ Now applying Theorem 4.5 yields the result.

\par{\bf Remark 7.2} In [S4] the author established the congruence
for $\sum_{k=0}^{p-1}\f{k^3\b{2k}k^3}{64^k}$ modulo $p^2$, where $p$
is a prime with $p>5$. See also [T3]. For the conjectures concerning
Theorems 7.4-7.6 see [S13, Conjectures 2.7, 2.11 and 2.16].

\section*{8. The congruence for
$\sum_{k=0}^{p-3}\b{a}k\b{-1-a}k\b{2k}{k}\f 1{4^k(k+2)}$ modulo
$p^3$} \pro{Lemma 8.1} For any positive integer $n$ and real number
$a\not=1,-2$ we have
$$\sum_{k=0}^{n-1}\b ak\b{-1-a}k\b{2k}k\f 1{4^k}\Big(\f 1{k+2}-\f 1{3(k+1)}+\f
1{3(a-1)(a+2)}\Big)=F(a,n)R(a,n),$$ where
$$\align &F(a,n)=\b an\b{-1-a}n\b{2n}n\f 1{4^n}\Big(\f 1{n+2}-\f 1{3(n+1)}+\f
1{3(a-1)(a+2)}\Big),
\\&R(a,n)=-\f{2n^2(n+2)}{n^2+(2a(a+1)-1)n+a(a+1)}.
\endalign$$
\endpro
Proof. It is easy to check that
$$R(a,0)=0\qtq{and} F(a,k)=F(a,k+1)R(a,k+1)-F(a,k)R(a,k).$$
Thus the result follows from (3.1).

\pro{Lemma 8.2} Let $p>3$ be a prime, $a\in\Bbb Z_p$, $a\not\e 0,\pm
1,-2\mod p$
 and $a'=(a-\ap)/p$. Then
 $$\f{\b a{p-2}\b{-1-a}{p-2}\b{2(p-2)}{p-2}}{4^{p-2}\cdot p}
 \e\f{2a'(a'+1)}{3a(a+1)(a-1)(a+2)}p^2\mod {p^3}.$$
 \endpro
 Proof. It is clear that
 $$\align &\b a{p-2}=\f{a(p-1)}{(a+1-p)(a+2-p)}\b{a-1}{p-1},\q
 \b{-1-a}{p-2}=\f{p-1}{1-a-p}\b{-1-a}{p-1},\\&
 \b{2(p-2)}{p-2}=\f{p(p-1)^2}{(2p-1)(2p-2)(2p-3)}\b{2p-1}{p-1}.\endalign$$
 Thus, applying [S9, Lemma 2.2] we get
$$\align&\f{\b a{p-2}\b{-1-a}{p-2}\b{2(p-2)}{p-2}}{4^{p-2}\cdot p}
\\&=\f{a(p-1)}{(a+1-p)(a+2-p)}\cdot\f{p-1}{1-a-p}\b{a-1}{p-1}\b{-1-a}{p-1}
\\&\q\times\f {4(p-1)^2}{(2p-1)(2p-2)(2p-3)\cdot 4^{p-1}}\b{2p-1}{p-1}
\\&\e \f a{(a+1)(a+2)(1-a)}\cdot\f{a'(a'+1)}{a^2}p^2\cdot\f 4{-6}
=\f{2a'(a'+1)}{3a(a+1)(a-1)(a+2)}p^2\mod {p^3}.\endalign$$ This
proves the lemma.

\pro{Theorem 8.1} Let $p>3$ be a prime, $a\in\Bbb Z_p$ and $a\not\e
0, \pm 1,-2\mod p$. Then
$$\align &\sum_{k=0}^{p-3}\b ak\b{-1-a}k\b{2k}k\f 1{4^k(k+2)}
\\&\e\f 13\sum_{k=0}^{p-2}\b ak\b{-1-a}k\b{2k}k\f 1{4^k(k+1)}-\f
1{3(a-1)(a+2)}S_p(a)
\\&\e
\f{a^2+a-3}{3(a-1)(a+2)}S_p(a)+\f{a+1}{3a}S_p(a+1) \mod
{p^3}.\endalign$$
\endpro
Proof. Set $a'=(a-\ap)/p$. Taking $n=p$ in Lemma 8.1 and noting that
$\b ap\b{-1-a}p\e -a'(a'+1)\mod p$ and $\b{2p}p\f 1{4^p}\e \f 12\mod
p$ yields
$$\align &\sum_{k=0}^{p-1}\b ak\b{-1-a}k\b{2k}k\f 1{4^k}\Big(\f 1{k+2}-\f 1{3(k+1)}+\f
1{3(a-1)(a+2)}\Big)
\\&\e \f{a'(a'+1)}{3(a-1)(a+2)}p^2\mod {p^3}.\endalign$$
Appealing to Lemmas 3.2, 8.2 and Theorem 3.1,
 $$\align&\sum_{k=0}^{p-3}\b
ak\b{-1-a}k\b{2k}k\f 1{4^k(k+2)}\\&\e \sum_{k=0}^{p-1}\b
ak\b{-1-a}k\b{2k}k\f 1{4^k(k+2)}-\f{2a'(a'+1)}{3a(a+1)(a-1)(a+2)}p^2
\\&\e \f 13\Big(\sum_{k=0}^{p-2}\b
ak\b{-1-a}k\b{2k}k\f 1{4^k(k+1)}-\f{a'(a'+1)}{a(a+1)}p^2\Big) -\f
1{3(a-1)(a+2)}S_p(a)
\\&\q+\f{a'(a'+1)}{3(a-1)(a+2)}p^2-\f{2a'(a'+1)}{3a(a+1)(a-1)(a+2)}p^2
\\&=\f 13\sum_{k=0}^{p-2}\b ak\b{-1-a}k\b{2k}k\f 1{4^k(k+1)}-\f
1{3(a-1)(a+2)}S_p(a)
\\&\e \f 13\Big(S_p(a)+\f{a+1}aS_p(a+1)\Big)-\f
1{3(a-1)(a+2)}S_p(a)\mod {p^3}.\endalign$$ This yields the result.

\pro{Theorem 8.2} Let $p>3$ be a prime. Then
$$\align \sum_{k=0}^{\f{p-1}2}\f{\b{2k}k^3}{64^k(k+2)}\e
\cases \f{52}{27}x^2-\f{26}{27}p-\f{p^2}{27x^2}\mod {p^3}\q\t{if
$p=x^2+4y^2\e 1\mod 4$,}
\\-\f{(p+1)^2}{3\cdot 2^{p-1}}
\b{(p-1)/2}{(p-3)/4}^2-p^2\Big(\f{13}{27}\b{(p-1)/2}{(p-3)/4}^{-2}
\\\qq+\f 16\b{(p-1)/2}{(p-3)/4}^2E_{p-3}\Big)\mod {p^3}\qq\t{if
$p\e 3\mod 4$.}
\endcases\endalign$$
\endpro
Proof. Since $p\mid \b{2k}k$ for $\f p2<k<p$, taking $a=-\f 12$ in
Theorem 8.1 gives
$$\align \sum_{k=0}^{(p-1)/2}\f{\b{2k}k^3}{64^k(k+2)}\e
\sum_{k=0}^{p-3}\f{\b{2k}k^3}{64^k(k+2)} \e \f
13\sum_{k=0}^{p-2}\f{\b{2k}k^3}{64^k(k+1)}+\f
4{27}\sum_{k=0}^{p-1}\f{\b{2k}k^3}{64^k}\mod {p^3}.\endalign$$
 Now applying (1.1) and Theorem 3.2 yields the result.

 \pro{Theorem 8.3} Let $p>5$ be a prime.
 For $p\e 1\mod 3$
and so $p=x^2+3y^2$ we have
$$\sum_{k=0}^{p-3}\f{\b{2k}k^2\b{3k}k}{108^k(k+2)}
\e \f{29}{15}x^2-\f{29}{30}p-\f{3p^2}{80x^2} \mod {p^3}.$$ For $p\e
2\mod 3$ we have
$$\align
&\sum_{k=0}^{p-3}\f{\b{2k}k^2\b{3k}k}{108^k(k+2)}
\\&\e -\f 23\b{(p-1)/2}{(p-5)/6}^2\Big(1+p\Big(2+\f 43\qp 2-\f 32\qp
3\Big) +p^2\Big(1+\f 83\qp 2+\f 29\qp 2^2\\&\q-3\qp 3-2\qp 2\qp
3+\f{15}8\qp 3^2+\f 34U_{p-3}\Big)\Big)-\f
{29}{120}p^2\b{(p-1)/2}{(p-5)/6}^{-2} \mod {p^3}.
\endalign$$
\endpro
 Proof. Set $a=-\f 13$. Then $3(a-1)(a+2)=-\f{20}3$.
 By Theorem 8.1,
$$\sum_{k=0}^{p-3}\f{\b{2k}k^2\b{3k}k}{108^k(k+2)}
\e \f 13\sum_{k=0}^{p-2}\f{\b{2k}k^2\b{3k}k}{108^k(k+1)} +\f
3{20}\sum_{k=0}^{p-1}\f{\b{2k}k^2\b{3k}k}{108^k} \mod {p^3}.$$
 Now applying (1.2) and Theorem 3.3 yields the
result.

\pro{Theorem 8.4} Let $p>7$ be a prime. Then
$$\sum_{k=0}^{p-3}\f{\b{2k}k^2\b{4k}{2k}}{256^k(k+2)}
\e\cases \f{17}{35}(4x^2-2p)-\f 4{105}p^2\mod {p^3}&\t{if
$p=x^2+2y^2\e 1,3\mod 8$,}
\\-\f 19R_2(p)\mod {p^2}&\t{if $p\e 5,7\mod 8$.}\endcases$$
\endpro
Proof. Set $a=-\f 14$. Then $3(a-1)(a+2)=-\f{105}{16}$.
 By Theorem 8.1,
$$\sum_{k=0}^{p-3}\f{\b{2k}k^2\b{4k}{2k}}{256^k(k+2)}\e\f 13
\sum_{k=0}^{p-2}\f{\b{2k}k^2\b{4k}{2k}}{256^k(k+1)}+\f{16}{105}
\sum_{k=0}^{p-1}\f{\b{2k}k^2\b{4k}{2k}}{256^k}\mod {p^3}.$$ Now,
applying (1.3) and Theorem 3.4 yields the result.

\pro{Theorem 8.5} Let $p$ be a prime with $p>11$. Then
$$\align &\sum_{k=0}^{p-3}\f{\b{2k}k\b{3k}k\b{6k}{3k}}{1728^k(k+2)}
\\&\e\cases\f{113}{231}
\sls p3(4x^2-2p)\mod {p^2}\qq\qq\qq\qq\q\ \t{if $p=x^2+4y^2\e 1\mod 4$,} \\
-\f 1{15}\b{\f{p-1}2}{[\f p{12}]}^2\big(1+p\big(\f 25-3\qp 2-\f
52\qp 3-\f 23H_{[\f p{12}]}\big)\big)\mod {p^2}\q\t{if $12\mid
p-7$,}
\\-\f 53\b{\f{p-1}2}{[\f p{12}]}^2\big(1+p\big(2-3\qp 2-\f
52\qp 3-\f 23H_{[\f p{12}]}\big)\big)\mod {p^2}\q\t{if $12\mid
p-11$.}\endcases\endalign$$
\endpro
Proof. Set $a=-\f 16$. Then $3(a-1)(a+2)=-\f{77}{12}$.
 By Theorem 8.1,
$$\sum_{k=0}^{p-3}\f{\b{2k}k\b{3k}k\b{6k}{3k}}{1728^k(k+2)}\e
\f 13\sum_{k=0}^{p-2}\f{\b{2k}k\b{3k}k\b{6k}{3k}}{1728^k(k+1)}
+\f{12}{77}\sum_{k=0}^{p-1}\f{\b{2k}k\b{3k}k\b{6k}{3k}}{1728^k}\mod
{p^3}.$$ Now applying Theorem 3.5 yields the result.

\section*{9. The congruence for
$\sum_{k=0}^{p-2}\b{a}k\b{-1-a}k\b{2k}{k}\f 1{4^k(k+1)^2}$ modulo
$p^3$} \pro{Lemma 9.1} For any positive integer $n$ and real number
$a\not=0$ we have
$$\align &\sum_{k=0}^{n-1}\b{a}k\b{-1-a}k\b{2k}k\f
1{4^k(k+1)^2}\\&=2S_n(a)+\f{2a^2+2a-1}{a^2}S_n(a+1)+F(a,n)R(a,n),\endalign$$
where
$$\align
&F(a,n)=\Big(\b{a}n\b{-1-a}n\big(\f 1{(n+1)^2}-2\big)
-\f{2a^2+2a-1}{a^2}\b{a+1}n\b{-2-a}n\Big)\b{2n}n\f 1{4^n},
\\&R(a,n)=-\f{2(n+1)^2((2a-1)(n+1)^2+(2-3a)(n+1)+a-1)}
{(2a-1)(n+1)^3+(4a^3+6a^2-a)(n+1)^2 +a^2(n+1)-a^2(a+2)}.\endalign$$
\endpro
Proof. It is easy to check that
$$R(a,0)=0\qtq{and} F(a,k)=F(a,k+1)R(a,k+1)-F(a,k)R(a,k).$$
Thus the result follows from (3.1).
 \pro{Theorem 9.1} Let $p>3$ be a
prime, $a\in\Bbb Z_p$, $a\not\e 0,-1,-2\mod p$ and $a'=(a-\ap)/p$.
Then
$$\align&\sum_{k=0}^{p-2}\b{a}k\b{-1-a}k\b{2k}k\f
1{4^k(k+1)^2}\\&\e
2S_p(a)+\f{2a^2+2a-1}{a^2}S_p(a+1)+\f{4a^3+6a^2-3a+2}{a^3(a+1)(a+2)}a'(a'+1)
p^3\mod {p^4}.\endalign$$
\endpro
Proof. Using [S9, Lemma 2.2] we see that
$$\b a{p-1}\b{-1-a}{p-1}=\f
a{a+1-p}\b{a-1}{p-1}\b{-a-1}{p-1} \e \f{a'(a'+1)}{a(a+1)}p^2\mod
{p^3}.\tag 9.1$$ Also, $\b{a+1}{p-1}\b{-2-a}{p-1}\e
\b{\ap+1}p\b{p-2-\ap}{p-1}\e 0\mod p$ and $\b{2(p-1)}{p-1}\f
1{4^{p-1}}=\f p{2p-1}\b{2p-1}{p-1}\f 1{4^{p-1}}\e -p\mod {p^2}.$
Thus, taking $n=p-1$ in Lemma 9.1 and then applying the above and
Lemma 3.2 we see that
$$\align &\sum_{k=0}^{p-2}\b{a}k\b{-1-a}k\b{2k}k\f
1{4^k(k+1)^2} \\&\e 2S_{p-1}(a)+\f{2a^2+2a-1}{a^2}S_{p-1}(a+1)
+\f{a'(a'+1)}{a(a+1)}\cdot(-p)\cdot\f{2(a-1)p^2}{a^2(a+2)}
\\&\e 2S_p(a)+2\f{a'(a'+1)}{a(a+1)}p^3+\f{2a^2+2a-1}{a^2}S_p(a+1)
\\&\q+\f{2a^2+2a-1}{a^2}\cdot\f{a'(a'+1)}{(a+1)(a+2)}p^3-\f{2(a-1)a'(a'+1)}
{a^3(a+1)(a+2)}p^3
\\&=2S_p(a)+\f{2a^2+2a-1}{a^2}S_p(a+1)+\f{4a^3+6a^2-3a+2}{a^3(a+1)(a+2)}a'(a'+1)
p^3\mod {p^4}.\endalign$$ This proves the theorem.

\pro{Theorem 9.2} Let $p>3$ be a prime, $a\in\Bbb Z_p$ and $a\not\e
0,-1,-2\mod p$. Then
$$\align&\sum_{k=0}^{p-2}\b{a}k\b{-1-a}k\b{2k}k\f
1{4^k(k+1)^2}\\&\e\Big(2-\f 1{a(a+1)}\Big)
\sum_{k=0}^{p-2}\b{a}k\b{-1-a}k\b{2k}k\f 1{4^k(k+1)}+\f
1{a(a+1)}S_p(a)\mod {p^3}.\endalign$$
\endpro
Proof. By Theorems 9.1 and 3.1,
$$\align&\sum_{k=0}^{p-2}\b{a}k\b{-1-a}k\b{2k}k\f
1{4^k(k+1)^2}\\&\e 2S_p(a)+\f{2a^2+2a-1}{a^2}S_p(a+1)\\&\e
2S_p(a)+\f{2a^2+2a-1}{a^2}\cdot\f
a{a+1}\Big(\sum_{k=0}^{p-2}\b{a}k\b{-1-a}k\b{2k}k\f
1{4^k(k+1)}-S_p(a)\Big)
\\&=\f{2a(a+1)-1}{a(a+1)}\sum_{k=0}^{p-2}\b{a}k\b{-1-a}k\b{2k}k\f 1{4^k(k+1)}+\f
1{a(a+1)}S_p(a)\mod {p^3}.\endalign$$ This proves the theorem.

\pro{Theorem 9.3} Let $p>3$ be a prime. Then
$$\sum_{k=0}^{\f{p-1}2}\f{\b{2k}k^3}{64^k(k+1)^2}\e
\cases 8x^2-4p+\f{p^2}{x^2}\mod {p^3}\qq\qq\t{if $p=x^2+4y^2\e 1\mod
4$,}
\\-\f{6(p+1)^2}{2^{p-1}}
\b{(p-1)/2}{(p-3)/4}^2-2p^2\b{(p-1)/2}{(p-3)/4}^{-2}
\\\qq-3p^2\b{(p-1)/2}{(p-3)/4}^2E_{p-3}\mod {p^3}\q\qq\qq\t{if
$4\mid p-3$}.
\endcases$$
\endpro
Proof. Since $p\mid \b{2k}k$ for $\f p2<k<p$, taking $a=-\f 12$ in
Theorem 9.2 gives
$$\sum_{k=0}^{\f{p-1}2}\f{\b{2k}k^3}{64^k(k+1)^2}\e
6\sum_{k=0}^{\f{p-1}2}\f{\b{2k}k^3}{64^k(k+1)}
-4\sum_{k=0}^{\f{p-1}2}\f{\b{2k}k^3}{64^k}\mod {p^3}.$$ Now applying
(1.1) and Theorem 3.2 yields the result.

\pro{Theorem 9.4} Let $p$ be a prime with $p>3$. For $p\e 1\mod 3$
and so $p=x^2+3y^2$ we have
$$\sum_{k=0}^{p-2}\f{\b{2k}k^2\b{3k}k}{108^k(k+1)^2}
\e 8x^2-4p+\f{9p^2}{8x^2} \mod {p^3}.$$ For $p\e 2\mod 3$ we have
$$\align
&\sum_{k=0}^{p-2}\f{\b{2k}k^2\b{3k}k}{108^k(k+1)^2}
\\&\e -13\b{(p-1)/2}{(p-5)/6}^2\Big(1+p\Big(2+\f 43\qp 2-\f 32\qp
3\Big) +p^2\Big(1+\f 83\qp 2+\f 29\qp 2^2\\&\q-3\qp 3-2\qp 2\qp
3+\f{15}8\qp 3^2+\f 34U_{p-3}\Big)\Big)-p^2\b{(p-1)/2}{(p-5)/6}^{-2}
\mod {p^3}.
\endalign$$
\endpro
Proof. Taking $a=-\f 13$ in Theorem 9.2 yields
$$\sum_{k=0}^{p-2}\f{\b{2k}k^2\b{3k}k}{108^k(k+1)^2}
\e \f{13}2\sum_{k=0}^{p-2}\f{\b{2k}k^2\b{3k}k}{108^k(k+1)} -\f
92\sum_{k=0}^{p-1}\f{\b{2k}k^2\b{3k}k}{108^k}\mod {p^3}.$$ Now
applying (1.2) and Theorem 3.3 yields the result.

\pro{Theorem 9.5} Let $p>3$ be a prime. Then
$$\sum_{k=0}^{p-2}\f{\b{2k}k^2\b{4k}{2k}}{256^k(k+1)^2}
\e\cases 8x^2-4p+\f{4p^2}{3x^2}\mod {p^3}&\t{if $p=x^2+2y^2\e
1,3\mod 8$,}
\\-\f {22}9R_2(p)\mod {p^2}&\t{if $p\e 5,7\mod 8$.}\endcases$$
\endpro
Proof. Taking $a=-\f 14$ in Theorem 9.2 yields
$$\sum_{k=0}^{p-2}\f{\b{2k}k^2\b{4k}{2k}}{256^k(k+1)^2}
\e\f{22}3\sum_{k=0}^{p-2}\f{\b{2k}k^2\b{4k}{2k}}{256^k(k+1)}
-\f{16}3\sum_{k=0}^{p-1}\f{\b{2k}k^2\b{4k}{2k}}{256^k}\mod {p^3}.$$
Now applying (1.3) and Theorem 3.4 yields the result.

\pro{Theorem 9.6} Let $p$ be a prime with $p>3$. Then
$$\align &\sum_{k=0}^{p-2}\f{\b{2k}k\b{3k}k\b{6k}{3k}}{1728^k(k+1)^2}
\\&\e\cases
\sls p3(8x^2-4p)\mod {p^2}\qq\qq\qq\qq\q\ \t{if $p=x^2+4y^2\e 1\mod 4$,} \\
-\f {46}{25}\b{\f{p-1}2}{[\f p{12}]}^2\big(1+p\big(\f 25-3\qp 2-\f
52\qp 3-\f 23H_{[\f p{12}]}\big)\big)\mod {p^2}\q\t{if $12\mid
p-7$,}
\\-46\b{\f{p-1}2}{[\f p{12}]}^2\big(1+p\big(2-3\qp 2-\f
52\qp 3-\f 23H_{[\f p{12}]}\big)\big)\mod {p^2}\q\t{if $12\mid
p-11$.}\endcases\endalign$$
\endpro
Proof. Taking $a=-\f 16$ in Theorem 9.2 yields
$$\sum_{k=0}^{p-2}\f{\b{2k}k\b{3k}k\b{6k}{3k}}{1728^k(k+1)^2}
\e
\f{46}{5}\sum_{k=0}^{p-2}\f{\b{2k}k\b{3k}k\b{6k}{3k}}{1728^k(k+1)}
-\f{36}5\sum_{k=0}^{p-1}\f{\b{2k}k\b{3k}k\b{6k}{3k}}{1728^k}\mod
{p^3}.$$ Now applying Theorem 3.5 yields the result.

\par{\bf Remark 9.1} Let $p>5$ be a prime. The congruence for
$\sum_{k=0}^{\f{p-1}2}\f{\b{2k}k^3}{64^k(k+1)^2}$ modulo $p^2$ can
be deduced from [T3,(4)]. In [S14] the author conjectured the
congruences modulo $p^2$ for the sums in Theorem 9.4-9.6.

\section*{10. The congruence for
$\sum_{k=0}^{p-2}\b{a}k\b{-1-a}k\b{2k}{k}\f 1{4^k(k+1)^3}$ modulo
$p^3$} \pro{Lemma 10.1} For any positive integer $n$ and real number
$a\not=0,-1$ we have
$$\align &\sum_{k=0}^{n-1}\b{a}k\b{-1-a}k\b{2k}k\f
1{4^k(k+1)^3}\\&=-\f
2{a(a+1)}+\f{(2a+1)^2}{a(a+1)}S_n(a)+\f{4a^2(a+1)^2-a(a+1)+1}
{a^3(a+1)}S_n(a+1)+F(a,n)R(a,n),\endalign$$ where
$$\align
&F(a,n)=\b{a}n\b{-1-a}n\Big(\f
1{(n+1)^3}-\f{(2a+1)^2}{a(a+1)}\Big)\b{2n}n\f 1{4^n}
\\&\qq\qq\q-\f{4a^2(a+1)^2-a(a+1)+1} {a^3(a+1)}\b{a+1}n\b{-2-a}n\b{2n}n\f
1{4^n},
\\&R(a,n)=-\f{2(n+1)^3((4a^2-2a+1)(n+1)^2+(-6a^2+3a-2)(n+1)+3a^2-a+1)}
{(4a^2-2a+1)(n+1)^4+(8a^4+12a^3+a)(n+1)^3+a^3(n+1)-a^3(a+2)}.
\endalign$$
\endpro
Proof. Using Maple it is easy to check that
$$F(a,k)=F(a,k+1)R(a,k+1)-F(a,k)R(a,k).$$
 On the other hand,
$$\align F(a,0)R(a,0)&=\Big(1-\f{(2a+1)^2}{a(a+1)}-\f{4a^2(a+1)^2-a(a+1)+1}
{a^3(a+1)}\Big)\\&\q\times\Big(-\f{2a^2}{7a^4+11a^3+4a^2-a+1}\Big)=\f
2{a(a+1)}.\endalign$$ Thus,
$$\align &\sum_{k=0}^{n-1}\b{a}k\b{-1-a}k\b{2k}k\f
1{4^k(k+1)^3}\\&\q-\f{(2a+1)^2}{a(a+1)}S_n(a)-\f{4a^2(a+1)^2-a(a+1)+1}
{a^3(a+1)}S_n(a+1)
\\&=\sum_{k=0}^{n-1}F(a,k)=\sum_{k=0}^{n-1}(F(a,k+1)R(a,k+1)-F(a,k)R(a,k))
\\&=F(a,n)R(a,n)-F(a,0)R(a,0)=F(a,n)R(a,n)-\f 2{a(a+1)}.\endalign$$
This proves the lemma.

\pro{Theorem 10.1} Let $p>3$ be a prime, $a\in\Bbb Z_p$ and $a\not\e
0,-1,-2\mod p$. Then
$$\align &\sum_{k=0}^{p-2}\b{a}k\b{-1-a}k\b{2k}k\f
1{4^k(k+1)^3}\\&\e -\f
2{a(a+1)}+\f{(2a+1)^2}{a(a+1)}S_{p}(a)+\f{4a^2(a+1)^2-a(a+1)+1}
{a^3(a+1)}S_{p}(a+1)
\\&\e -\f
2{a(a+1)}+\f{4a^2(a+1)^2-a(a+1)+1}
{a^2(a+1)^2}\sum_{k=0}^{p-2}\b{a}k\b{-1-a}k\b{2k}k\f 1{4^k(k+1)}
\\&\qq+\f{2a(a+1)-1}{a^2(a+1)^2}S_p(a)\mod{p^3}.\endalign$$
\endpro

Proof. Clearly $\b a{p-1}\b{-1-a}{p-1}\e 0\mod {p^2}$ and $p\mid
\b{2(p-1)}{p-1}$. Hence, taking $n=p-1$ in Lemma 10.1 yields
$$\align &\sum_{k=0}^{p-2}\b{a}k\b{-1-a}k\b{2k}k\f
1{4^k(k+1)^3}\\&\e-\f
2{a(a+1)}+\f{(2a+1)^2}{a(a+1)}S_{p-1}(a)+\f{4a^2(a+1)^2-a(a+1)+1}
{a^3(a+1)}S_{p-1}(a+1).\endalign$$ By Lemma 3.2, $S_{p-1}(a)\e
S_{p}(a)\mod {p^3}$ and $S_{p-1}(a+1)\e S_{p}(a+1)\mod {p^3}$. Now,
from the above and Theorem 3.1 we get $$\align
&\sum_{k=0}^{p-2}\b{a}k\b{-1-a}k\b{2k}k\f 1{4^k(k+1)^3}\\&\e-\f
2{a(a+1)}+\f{(2a+1)^2}{a(a+1)}S_{p}(a)+\f{4a^2(a+1)^2-a(a+1)+1}
{a^3(a+1)}\\&\q\times\f
a{a+1}\Big(\sum_{k=0}^{p-2}\b{a}k\b{-1-a}k\b{2k}k\f
1{4^k(k+1)}-S_p(a)\Big)\mod {p^3}.\endalign$$ This yields the
result.

\pro{Theorem 10.2} Let $p$ be a prime with $p>3$. For $p\e 1\mod 3$
and so $p=x^2+3y^2$ we have
$$\sum_{k=0}^{p-2}\f{\b{2k}k^2\b{3k}k}{108^k(k+1)^3}
\e 9-2x^2+p+\f{117p^2}{16x^2} \mod {p^3}.$$ For $p\e 2\mod 3$ we
have
$$\align
&\sum_{k=0}^{p-2}\f{\b{2k}k^2\b{3k}k}{108^k(k+1)^3}
\\&\e 9-\f{115}2\b{(p-1)/2}{(p-5)/6}^2\Big(1+p\Big(2+\f 43\qp 2-\f 32\qp
3\Big) +p^2\Big(1+\f 83\qp 2+\f 29\qp 2^2\\&\q-3\qp 3-2\qp 2\qp
3+\f{15}8\qp 3^2+\f 34U_{p-3}\Big)\Big)+\f
14p^2\b{(p-1)/2}{(p-5)/6}^{-2} \mod {p^3}.
\endalign$$
\endpro
Proof. Putting $a=-\f 13$ in Theorem 10.1 gives
$$\sum_{k=0}^{p-2}\f{\b{2k}k^2\b{3k}k}{108^k(k+1)^3}
\e 9+\f{115}4\sum_{k=0}^{p-2}\f{\b{2k}k^2\b{3k}k}{108^k(k+1)}
-\f{117}4\sum_{k=0}^{p-1}\f{\b{2k}k^2\b{3k}k}{108^k}\mod {p^3}.$$
Now applying (1.2) and Theorem 3.3 yields the result.

\pro{Theorem 10.3} Let $p>3$ be a prime. Then
$$\sum_{k=0}^{p-2}\f{\b{2k}k^2\b{4k}{2k}}{256^k(k+1)^3}
\e\cases \f{32}3-\f{16}3x^2+\f 83p+\f{88p^2}{9x^2}\mod {p^3}&\t{if
$p=x^2+2y^2\e 1,3\mod 8$,}
\\\f{32}3-\f {340}{27}R_2(p)\mod {p^2}&\t{if $p\e 5,7\mod 8$.}\endcases$$
\endpro
Proof. Putting $a=-\f 14$ in Theorem 10.1 gives
$$\sum_{k=0}^{p-2}\f{\b{2k}k^2\b{4k}{2k}}{256^k(k+1)^3}
\e \f{32}3+\f{340}9\sum_{k=0}^{p-2}
\f{\b{2k}k^2\b{4k}{2k}}{256^k(k+1)}
-\f{352}9\sum_{k=0}^{p-1}\f{\b{2k}k^2\b{4k}{2k}}{256^k} \mod
{p^3}.$$ Now applying (1.3) and Theorem 3.4 yields the result.

\pro{Theorem 10.4} Let $p$ be a prime with $p>5$. Then
$$\align &\sum_{k=0}^{p-2}\f{\b{2k}k\b{3k}k\b{6k}{3k}}{1728^k(k+1)^3}
\\&\e\cases\f{72}5-\f{16}5
\sls p3(4x^2-2p)\mod {p^2}\qq\qq\qq\qq\q\ \t{if $p=x^2+4y^2\e 1\mod 4$,} \\
\f{72}5-\f {1576}{125}\b{\f{p-1}2}{[\f p{12}]}^2\big(1+p\big(\f
25-3\qp 2-\f 52\qp 3-\f 23H_{[\f p{12}]}\big)\big)\mod {p^2}\q\t{if
$12\mid p-7$,}
\\\f{72}5-\f{1576}5\b{\f{p-1}2}{[\f p{12}]}^2\big(1+p\big(2-3\qp 2-\f
52\qp 3-\f 23H_{[\f p{12}]}\big)\big)\mod {p^2}\q\t{if $12\mid
p-11$.}\endcases\endalign$$
\endpro
Proof. Putting $a=-\f 16$ in Theorem 10.1 gives
$$\sum_{k=0}^{p-2}\f{\b{2k}k\b{3k}k\b{6k}{3k}}{1728^k(k+1)^3}
\e
\f{72}5+\f{1576}{25}\sum_{k=0}^{p-2}\f{\b{2k}k\b{3k}k\b{6k}{3k}}{1728^k(k+1)}
-\f{1656}{25}\sum_{k=0}^{p-1}\f{\b{2k}k\b{3k}k\b{6k}{3k}}{1728^k}
\mod {p^3}.$$ Now applying Theorem 3.5 yields the result.

\pro{Theorem 10.5} Let $p$ be an odd prime. Then
$$\align &\sum_{k=0}^{p-1}\f{\b{2k}k^3}{64^k(2k-1)^3}
\\&\e\cases \f{3p^2}{4x^2}\mod
{p^3}&\t{if $p=x^2+4y^2\e 1\mod 4$,}
\\-\f{3(p+1)^2}{2^{p-1}}\b{\f{p-1}2}{\f{p-3}4}^2
-\f 32p^2\b{\f{p-1}2}{\f{p-3}4}^2E_{p-3}\mod {p^3}&\t{if $p\e 3\mod
4$.}\endcases\endalign$$
 \endpro
 Proof. Since $\f 1{2k-1}\b{2k}k=\f 2k\b{2(k-1)}{k-1}$, we see that
 $$\sum_{k=1}^{p-1}\f{\b{2k}k^3}{m^k(2k-1)^3}=\sum_{k=1}^{p-1}
 \f{8\b{2(k-1)}{k-1}^3}{m^kk^3}=\f
 8m\sum_{r=0}^{p-2}\f{\b{2r}r^3}{m^r(r+1)^3}.\tag 10.1$$
 Hence, appealing to Theorem 10.1 (with $a=-\f 12$) we obtain
 $$\align \sum_{k=0}^{p-1}\f{\b{2k}k^3}{64^k(2k-1)^3}
 &=-1+\f 8{64}\sum_{k=0}^{p-2}\f{\b{2k}k^3}{64^k(k+1)^3}
\\&\e -1+\f 18\Big(8+24\sum_{k=0}^{p-2}\f{\b{2k}k^3}{64^k(k+1)}-24
\sum_{k=0}^{p-1}\f{\b{2k}k^3}{64^k}\Big)
\\&=3\sum_{k=0}^{p-2}\f{\b{2k}k^3}{64^k(k+1)}
-3\sum_{k=0}^{p-1}\f{\b{2k}k^3}{64^k}\mod {p^3}.\endalign$$ Now
applying (1.1) and Theorem 3.2 yields the result.

\section*{11. The congruence for
$\sum_{k=0}^{p-4}\b{a}k\b{-1-a}k\b{2k}{k}\f 1{4^k(k+3)}$ modulo
$p^3$} \pro{Lemma 11.1} For any positive integer $n$ and real number
$a\not=0,\pm 1,\pm 2,-3$ we have
$$\align &\sum_{k=0}^{n-1}\b{a}k\b{-1-a}k\b{2k}{k}\f 1{4^k(k+3)}
\\&=\f{3a(a+1)-10}{15(a-1)(a+2)}S_n(a)+\f{(a+1)(3a(a+1)-16)}{15a(a-2)(a+3)}S_n(a+1)
+G(a,n),
\endalign$$
where
$$G(a,n)=\f{2n^2N(a,n)}
{15a(a^2-1)(a^2-4)(a+3)(n+1)(n+2)}\b{a+1}n\b{-1-a}n\b{2n}n\f
1{4^n}$$ and
$$N(a,n)=(3a^4+12a^3-13a^2-50a+32)n^2+(3a^4+12a^3-29a^2-82a+96)n+64.$$
\endpro
Proof. Set
$$\align&F(a,n)=\b an\b{-1-a}n\b{2n}n\f 1{4^n}\Big(\f 1{n+3}
-\f{3a(a+1)-10}{15(a-1)(a+2)}
\\&\qq\qq\q-\f{(a+1)(3a(a+1)-16)}{15a(a-2)(a+3)}
\cdot\f{a+1+n}{a+1-n}\Big).\endalign$$ Using Maple it is easy to
check that $F(a,k)=G(a,k+1)-G(a,k).$ Note that
$\b{a+1}n\b{-2-a}n=\f{a+1+n}{a+1-n}\b an\b{-1-a}n$. We then have
$$\align &\sum_{k=0}^{n-1}\b{a}k\b{-1-a}k\b{2k}{k}\f 1{4^k(k+3)}
\\&\q-\Big(\f{3a(a+1)-10}{15(a-1)(a+2)}S_n(a)+\f{(a+1)(3a(a+1)-16)}
{15a(a-2)(a+3)}S_n(a+1)\Big)
\\&=\sum_{k=0}^{n-1}F(a,k)=\sum_{k=0}^{n-1}(G(a,k+1)-G(a,k))=G(a,n)-G(a,0)
=G(a,n).\endalign$$ This proves the lemma.

 \pro{Lemma 11.2} Let $p>5$ be a prime, $a\in\Bbb Z_p$,
 $a\not\e 0,\pm 1,\pm 2,-3\mod p$ and $a'=(a-\ap)/p$. Then
$$\f{\b a{p-3}\b{-1-a}{p-3}\b{2(p-3)}{p-3}}{4^{p-3}\cdot p}
\e-\f{32a'(a'+1)p^2}{15a(a^2-1)(a^2-4)(a+3)}\mod {p^3}.$$
\endpro
Proof. It is easy to see that
$$\align &\b a{p-3}\b{-1-a}{p-3}=\f{(p-2)^2}{(a+3-p)(2-a-p)}\b a{p-2}\b{-1-a}{p-2}
\\&\b{2(p-3)}{p-3}=\f{(p-2)^2}{(2p-4)(2p-5)}\b{2(p-2)}{p-2}.\endalign$$
Thus, appealing to Lemma 8.2 we obtain
$$\align \f{\b a{p-3}\b{-1-a}{p-3}\b{2(p-3)}{p-3}}{4^{p-3}\cdot p}
&=\f {4(p-2)^4}{(a+3-p)(2-a-p)(2p-4)(2p-5)}\cdot \f{\b
a{p-2}\b{-1-a}{p-2}\b{2(p-2)}{p-2}}{4^{p-2}\cdot p}
\\&\e \f{16}{5(a+3)(2-a)}\cdot\f{2a'(a'+1)}{3a(a+1)(a-1)(a+2)}p^2
\\&=-\f{32a'(a'+1)p^2}{15(a-2)(a-1)a(a+1)(a+2)(a+3)}\mod {p^3}.
\endalign$$
This proves the lemma.

 \pro{Theorem 11.1} Let $p>5$ be a prime and $a\in\Bbb Z_p$ with
 $a\not\e 0,\pm 1,\pm 2,-3\mod p$. Then
 $$\align &\sum_{k=0}^{p-4}\b{a}k\b{-1-a}k\b{2k}{k}\f 1{4^k(k+3)}
 \\&\e \f{3a(a+1)-10}{15(a-1)(a+2)}S_p(a)+\f{(a+1)(3a(a+1)-16)}{15a(a-2)(a+3)}
 S_p(a+1)\\&\e \f{28-6a(a+1)}{15(a-1)(a^2-4)(a+3)}S_p(a)
\\&\qq+\f{3a(a+1)-16}{15(a-2)(a+3)}\sum_{k=0}^{p-2}\b{a}k\b{-1-a}k\b{2k}{k}\f
1{4^k(k+1)}\mod {p^3}.\endalign$$
 \endpro
Proof. Set $a'=(a-\ap)/p$. From Lemmas 3.2, 8.2 and 11.2 we deduce
that
$$\align&\sum_{k=0}^{p-1}\b{a}k\b{-1-a}k\b{2k}{k}\f 1{4^k(k+3)}
\\&\e \sum_{k=0}^{p-4}\b{a}k\b{-1-a}k\b{2k}{k}\f 1{4^k(k+3)}
-\f{32a'(a'+1)p^2}{15a(a^2-1)(a^2-4)(a+3)}\mod {p^3}.\endalign$$ Let
$G(a,n)$ be given in Lemma 11.1. Then
$$G(a,p)\e \f{64p^2}{15a(a^2-1)(a^2-4)(a+3)}
\b{a+1}p\b{-1-a}p\b{2p}p\f 1{4^p}\mod {p^3}.$$ Since $\b{2p}p\f
1{4^p}=\f 24\b{2p-1}{p-1}\f 1{4^{p-1}}\e \f 12\mod p$ and
$$\b{a+1}p\b{-1-a}p=\b{\ap+1+a'p}p\b{p-1-\ap-(a'+1)p}p\e
a'(-a'-1)\mod p,$$ we see that
$$G(a,p)\e -\f{32a'(a'+1)p^2}{15a(a^2-1)(a^2-4)(a+3)}\mod {p^3}.$$
Now, taking $n=p$ in Lemma 11.1 and then applying the above and
Theorem 3.1 we deduce that
$$\align &\sum_{k=0}^{p-4}\b{a}k\b{-1-a}k\b{2k}{k}\f 1{4^k(k+3)}
\\&\e \sum_{k=0}^{p-1}\b{a}k\b{-1-a}k\b{2k}{k}\f 1{4^k(k+3)}
+\f{32a'(a'+1)p^2}{15a(a^2-1)(a^2-4)(a+3)}
\\&=\f{3a(a+1)-10}{15(a-1)(a+2)}S_p(a)+\f{(a+1)(3a(a+1)-16)}{15a(a-2)(a+3)}
 S_p(a+1)\\&\qq+G(a,p)+\f{32a'(a'+1)p^2}{15a(a^2-1)(a^2-4)(a+3)}
 \\&\e \f{3a(a+1)-10}{15(a-1)(a+2)}S_p(a)+\f{(a+1)(3a(a+1)-16)}{15a(a-2)(a+3)}
 S_p(a+1)
\\&\e  \f{3a(a+1)-10}{15(a-1)(a+2)}S_p(a)
\\&\qq+\f{3a(a+1)-16}{15(a-2)(a+3)}
 \Big(\sum_{k=0}^{p-2}\b{a}k\b{-1-a}k\b{2k}{k}\f
1{4^k(k+1)}-S_p(a)\Big)
 \mod {p^3},\endalign$$
 which yields the result.

 \pro{Theorem 11.2} Let $p>5$ be a prime. Then
$$\sum_{k=0}^{\f{p-1}2}\f{\b{2k}k^3}{64^k(k+3)}\e
\cases \f{172}{135}x^2-\f{86}{135}p-\f{118p^2}{3375x^2}\mod
{p^3}\q\t{if $p=x^2+4y^2\e 1\mod 4$,}
\\-\f{67(p+1)^2}{375\cdot 2^{p-1}}
\b{(p-1)/2}{(p-3)/4}^2-\f{43}{135}p^2\b{(p-1)/2}{(p-3)/4}^{-2}
\\\qq-\f {67}{750}p^2\b{(p-1)/2}{(p-3)/4}^2E_{p-3}\mod {p^3}\qq\q\qq\t{if
$4\mid p-3$}.
\endcases$$
\endpro
Proof. Since $p\mid \b{2k}k$ for $\f p2<k<p$, taking $a=-\f 12$ in
Theorem 11.1 gives
$$\sum_{k=0}^{\f{p-1}2}\f{\b{2k}k^3}{64^k(k+3)}
\e \sum_{k=0}^{p-4}\f{\b{2k}k^3}{64^k(k+3)} \e \f{472}{3375}
\sum_{k=0}^{p-1}\f{\b{2k}k^3}{64^k}+\f{67}{375}\sum_{k=0}^{(p-1)/2}
\f{\b{2k}k^3}{64^k(k+1)}\mod {p^3}.$$ Now applying (1.1) and Theorem
3.2 yields the result.

\pro{Theorem 11.3} Let $p$ be a prime with $p>7$. For $p\e 1\mod 3$
and so $p=x^2+3y^2$ we have
$$\sum_{k=0}^{p-4}\f{\b{2k}k^2\b{3k}k}{108^k(k+3)}
\e \f{32}{25}x^2-\f{16}{25}p-\f{99p^2}{2800x^2} \mod {p^3}.$$ For
$p\e 2\mod 3$ we have
$$\align
&\sum_{k=0}^{p-4}\f{\b{2k}k^2\b{3k}k}{108^k(k+3)}
\\&\e -\f 5{14}\b{(p-1)/2}{(p-5)/6}^2\Big(1+p\Big(2+\f 43\qp 2-\f 32\qp
3\Big) +p^2\Big(1+\f 83\qp 2+\f 29\qp 2^2\\&\q-3\qp 3-2\qp 2\qp
3+\f{15}8\qp 3^2+\f 34U_{p-3}\Big)\Big)-\f
4{25}p^2\b{(p-1)/2}{(p-5)/6}^{-2} \mod {p^3}.
\endalign$$
\endpro
Proof. Taking $a=-\f 13$ in Theorem 11.1 gives
$$\sum_{k=0}^{p-4}\f{\b{2k}k^2\b{3k}k}{108^k(k+3)}
\e \f{99}{700}\sum_{k=0}^{p-1}\f{\b{2k}k^2\b{3k}k}{108^k} +\f 5{28}
\sum_{k=0}^{p-2}\f{\b{2k}k^2\b{3k}k}{108^k(k+1)}\mod {p^3}.$$ Now
applying (1.2) and Theorem 3.3 deduces the result.

\pro{Theorem 11.4} Let $p>11$ be a prime. Then
$$\sum_{k=0}^{p-4}\f{\b{2k}k^2\b{4k}{2k}}{256^k(k+3)}
\e\cases \f{676}{525}x^2-\f{338}{525}p-\f{1864p^2}{51975x^2}\mod
{p^3}&\t{if $p=x^2+2y^2\e 1,3\mod 8$,}
\\-\f {53}{891}R_2(p)\mod {p^2}&\t{if $p\e 5,7\mod 8$.}\endcases$$
\endpro
Proof. Taking $a=-\f 14$ in Theorem 11.1 gives
$$\sum_{k=0}^{p-4}\f{\b{2k}k^2\b{4k}{2k}}{256^k(k+3)}
\e \f{7456}{51975}\sum_{k=0}^{p-1}\f{\b{2k}k^2\b{4k}{2k}}{256^k}
+\f{53}{297} \sum_{k=0}^{p-2}\f{\b{2k}k^2\b{4k}{2k}}{256^k(k+1)}
 \mod {p^3}.$$
Now applying (1.3) and Theorem 3.4 deduces the result.

\pro{Theorem 11.5} Let $p$ be a prime with $p>17$. Then
$$\align &\sum_{k=0}^{p-4}\f{\b{2k}k\b{3k}k\b{6k}{3k}}{1728^k(k+3)}
\\&\e\cases
\f{25}{77}\sls p3(4x^2-2p)\mod {p^2}\qq\qq\qq\qq\q\ \t{if $p=x^2+4y^2\e 1\mod 4$,} \\
-\f {197}{5525}\b{\f{p-1}2}{[\f p{12}]}^2\big(1+p\big(\f 25-3\qp
2-\f 52\qp 3-\f 23H_{[\f p{12}]}\big)\big)\mod {p^2}\q\t{if $12\mid
p-7$,}
\\-\f{197}{221}\b{\f{p-1}2}{[\f p{12}]}^2\big(1+p\big(2-3\qp 2-\f
52\qp 3-\f 23H_{[\f p{12}]}\big)\big)\mod {p^2}\q\t{if $12\mid
p-11$.}\endcases\endalign$$
\endpro
Proof. Taking $a=-\f 16$ in Theorem 11.1 gives
$$\sum_{k=0}^{p-4}\f{\b{2k}k\b{3k}k\b{6k}{3k}}{1728^k(k+3)}
\e
\f{12456}{85085}\sum_{k=0}^{p-1}\f{\b{2k}k\b{3k}k\b{6k}{3k}}{1728^k}
+\f{197}{1105}\sum_{k=0}^{p-2}\f{\b{2k}k\b{3k}k\b{6k}{3k}}{1728^k(k+1)}\mod
{p^3}.$$
 Now applying Theorem 3.5 deduces the result.

\section*{12. The congruence for
$\sum_{k=0}^{p-1}\b{a}k\b{-1-a}k\b{2k}{k}\f 1{4^k(a-1+k)}$ modulo
$p^3$} \pro{Lemma 12.1} For any positive integer $n$ and real number
$a\not=0,1$ we have
$$\align &\sum_{k=0}^{n-1}\b{a}k\b{-1-a}k\b{2k}k\f
1{4^k(a-1+k)}\\&=\f{(a-1)^2+1}{2(a-1)^3}S_n(a)+\f 1
{2(a-1)}S_n(a-1)+F(a,n)R(a,n),\endalign$$ where
$$\align
&F(a,n)=\Big(\b{a}n\b{-1-a}n\big(\f
1{a-1+n}-\f{(a-1)^2+1}{2(a-1)^3}\big)\\&\qq\qq-\f
1{2(a-1)}\b{a-1}n\b{-a}n\Big)\b{2n}n\f 1{4^n},
\\&R(a,n)=-\f{2n^3}{n^2+(2a^2-2a+1)n+a(a-1)}.\endalign$$
\endpro
Proof. It is easy to check that
$$R(a,0)=0\qtq{and} F(a,k)=F(a,k+1)R(a,k+1)-F(a,k)R(a,k).$$
Thus the result follows from (3.1).
 \pro{Theorem 12.1} Let $p$ be an
odd prime and $a\in\Bbb Z_p$ with $a\not\e 0,\pm 1\mod p$. Then
$$\align &\sum_{k=0}^{p-1}\b ak\b{-1-a}k\b{2k}k\f 1{4^k(a+k)}
\e \f 1{2a}S_p(a)+\f{(a+1)^2}{2a^3}S_p(a+1)\\&\e\f{a+1}{2a^2}
\sum_{k=0}^{p-2}\b ak\b{-1-a}k\b{2k}k\f 1{4^k(k+1)}-\f 1{2a^2}S_p(a)
\mod {p^3},
\\&\sum_{k=0}^{p-1}\b ak\b{-1-a}k\b{2k}k\f 1{4^k(a+k-1)}
\e\f{(a-1)^2+1}{2(a-1)^3}S_p(a)+\f{(a+1)^2}{2a^2(a-1)}S_p(a+1)
\\&\e \f{a+1}{2a(a-1)}\sum_{k=0}^{p-2}\b ak\b{-1-a}k\b{2k}k\f 1{4^k(k+1)}
-\f{a^2-3a+1}{2a(a-1)^3}S_p(a) \mod {p^3}.\endalign$$
\endpro
Proof. Since
$$\b{a-1}k\b{-a}k=\f {a-k}a\b ak\f
a{a+k}\b{-1-a}k=\Big(\f{2a}{a+k}-1\Big)\b ak\b{-1-a}k,$$ using Lemma
2.1 and Theorem 3.1 we see that
$$\align &2a\sum_{k=0}^{p-1}\b ak\b{-1-a}k\b{2k}k\f 1{4^k(a+k)}
\\&=S_p(a)+S_p(a-1)\e S_p(a)+\f{(a+1)^2}{a^2}S_p(a+1)
\\&\e\f{a+1}{a} \sum_{k=0}^{p-2}\b ak\b{-1-a}k\b{2k}k\f
1{4^k(k+1)}-\f 1{a}S_p(a)\mod {p^3}.\endalign$$ On the other hand,
taking $n=p$ in Lemma 12.1 and then applying Theorem 3.1 yields the
remaining part.

\pro{Theorem 12.2} Let $p>3$ be a prime. Then
$$\sum_{k=0}^{p-1}\f{\b{2k}k^3}{64^k(2k-3)}\e
\cases -\f{26}{27}x^2+\f{13}{27}p+\f{11p^2}{108x^2}\mod {p^3}\q\t{if
$p=x^2+4y^2\e 1\mod 4$,}
\\-\f{(p+1)^2}{6\cdot 2^{p-1}}
\b{(p-1)/2}{(p-3)/4}^2+\f{13}{54}p^2\b{(p-1)/2}{(p-3)/4}^{-2}
\\\qq-\f 1{12}p^2\b{(p-1)/2}{(p-3)/4}^2E_{p-3}\mod {p^3}\q\qq\qq\t{if
$4\mid p-3$}.
\endcases$$
\endpro
Proof. Taking $a=-\f 12$ in Theorem 12.1 yields
$$\sum_{k=0}^{p-1}\f{\b{2k}k^3}{64^k(k-\f 32)}
\e \f 13\sum_{k=0}^{p-2}\f{\b{2k}k^3}{64^k(k+1)}-\f{22}{27}
\sum_{k=0}^{p-1}\f{\b{2k}k^3}{64^k}\mod {p^3}.$$ This together with
(1.1) and Theorem 3.2 yields the result.

\section*{13. Some challenging conjectures}
\par Calculations by Maple suggests the following conjectures:
\pro{Conjecture 13.1} Let $p$ be a prime with $p>7$.
\par $(\t{\rm i})$ If $p\e 1,2,4\mod 7$ and so $p=x^2+7y^2$, then
$$\align &\sum_{k=0}^{(p-1)/2}\f{\b{2k}k^3}{k+2}
 \e -\f{466}{27}y^2+\f{29}{27}p+\f{17p^2}{864y^2}\mod {p^3},
\\&\sum_{k=0}^{(p-1)/2}\f{\b{2k}k^3}{k+3}
 \e -\f{36052}{3375}y^2+\f{2378}{3375}p+\f{1421p^2}{108000y^2}\mod
 {p^3},
\\&\sum_{k=0}^{(p-1)/2}\f{\b{2k}k^3}{(k+1)^2}\e
-68y^2+p-\f{p^2}{4y^2}\mod {p^3},
\\&\sum_{k=0}^{(p-1)/2}\f{\b{2k}k^3}{(k+2)^2}\e
-\f{287}{27}y^2+\f{119}{216}p+\f{7p^2}{864y^2}\mod {p^3},
\\&\sum_{k=0}^{(p-1)/2}\f{\b{2k}k^3}{(k+3)^2}\e
-\f{68722}{16875}y^2+\f{33079}{135000}p+\f{4987p^2}{1080000y^2}\mod
{p^3},
\\&\sum_{k=0}^{(p-1)/2}\f{\b{2k}k^3}{(k+1)^3}\e
\f 18-\f{201}2y^2-\f 94p-\f{39p^2}{32y^2}\mod {p^3}.
\endalign$$
\par $(\t{\rm ii})$ If $p\e 3,5,6\mod 7$, then
$$\align &\sum_{k=0}^{(p-1)/2}\f{\b{2k}k^3}{k+2}
\e \f 1{24}R_7(p)-\f p9\mod {p^2},
\\&\sum_{k=0}^{(p-1)/2}\f{\b{2k}k^3}{k+3}
 \e \f 1{120}R_7(p)-\f {11}{225}p\mod {p^2},
\\&\sum_{k=0}^{(p-1)/2}\f{\b{2k}k^3}{(k+2)^2}
\e \f 7{72}R_7(p)-\f 7{72}p\mod {p^2},
\\&\sum_{k=0}^{(p-1)/2}\f{\b{2k}k^3}{(k+3)^2}
\e \f {37}{3600}R_7(p)-\f {307}{9000}p\mod {p^2},
\\&\sum_{k=0}^{(p-1)/2}\f{\b{2k}k^3}{(k+1)^3}
\e\f 18+24R_7(p)+18p\mod {p^2}.
 \endalign$$
\endpro

\pro{Conjecture 13.2} Let $p$ be a prime with $p>5$.
\par $(\t{\rm i})$ If $p\e 1\mod 3$ and so $p=x^2+3y^2$, then
$$\align
&\sum_{k=0}^{(p-1)/2}\f{\b{2k}k^3}{16^k(k+2)}
 \e \f{64}{27}x^2-\f 43p-\f{p^2}{9x^2}\mod {p^3},
 \\&\sum_{k=0}^{(p-1)/2}\f{\b{2k}k^3}{16^k(k+3)}
 \e  \f{5056}{3375}x^2-\f {302}{375}p-\f{97p^2}{1125x^2}\mod {p^3},
\\&\sum_{k=0}^{(p-1)/2}\f{\b{2k}k^3}{16^k(k+1)^2}\e
-24y^2+2p-\f{p^2}{2y^2}\mod {p^3},
\\&\sum_{k=0}^{(p-1)/2}\f{\b{2k}k^3}{16^k(k+2)^2}\e
-\f{40}9y^2+\f{14}{27}p-\f{p^2}{54y^2}\mod {p^3},
\\&\sum_{k=0}^{(p-1)/2}\f{\b{2k}k^3}{16^k(k+3)^2}\e
-\f{9656}{5625}y^2+\f{3958}{16875}p+\f{169p^2}{33750y^2}\mod {p^3},
\\&\sum_{k=0}^{(p-1)/2}\f{\b{2k}k^3}{16^k(k+1)^3}\e
2-24y^2-\f{5p^2}{2y^2}\mod {p^3}.
\endalign$$
\par $(\t{\rm ii})$ If $p\e 2\mod 3$, then
$$\align
&\sum_{k=0}^{(p-1)/2}\f{\b{2k}k^3}{16^k(k+2)}
 \e -\f 4{27}R_3(p)-\f 4{27}p\mod {p^2},
 \\&\sum_{k=0}^{(p-1)/2}\f{\b{2k}k^3}{16^k(k+3)}
 \e -\f 4{135}R_3(p)-\f {38}{675}p\mod {p^2},
\\&\sum_{k=0}^{(p-1)/2}\f{\b{2k}k^3}{16^k(k+2)^2}\e
-\f{16}{27}R_3(p)-\f 29p\mod {p^2},
\\&\sum_{k=0}^{(p-1)/2}\f{\b{2k}k^3}{16^k(k+3)^2}\e
-\f{56}{675}R_3(p)-\f{58}{1125}p\mod {p^2},
\\&\sum_{k=0}^{(p-1)/2}\f{\b{2k}k^3}{16^k(k+1)^3}\e
2-32R_3(p)-4p\mod {p^2}.
\endalign$$
\endpro

\pro{Conjecture 13.3} Let $p$ be a prime with $p>5$.
\par $(\t{\rm i})$ If $p\e 1\mod 4$ and so $p=x^2+4y^2$, then
$$\align
&\sum_{k=0}^{(p-1)/2}\f{\b{2k}k^3}{(-8)^k(k+2)}
 \e \f{64}{27}x^2-\f{35}{27}p-\f{4p^2}{27x^2}\mod {p^3},
 \\&\sum_{k=0}^{(p-1)/2}\f{\b{2k}k^3}{(-8)^k(k+3)}
 \e \f{202}{135}x^2-\f{538}{675}p-\f{611p^2}{6750x^2}\mod {p^3},
 \\&\sum_{k=0}^{(p-1)/2}\f{\b{2k}k^3}{(-8)^k(k+1)^2}\e
-32y^2+p-\f{7p^2}{16y^2}\mod {p^3},
\\&\sum_{k=0}^{(p-1)/2}\f{\b{2k}k^3}{(-8)^k(k+2)^2}\e
-\f{152}{27}y^2+\f{16}{27}p+\f{7p^2}{216y^2}\mod {p^3},
\\&\sum_{k=0}^{(p-1)/2}\f{\b{2k}k^3}{(-8)^k(k+3)^2}\e
-\f{1504}{675}y^2+\f{811}{3375}p+\f{1927p^2}{270000y^2}\mod {p^3},
\\&\sum_{k=0}^{(p-1)/2}\f{\b{2k}k^3}{(-8)^k(k+1)^3}\e
-1-48y^2-\f{33p^2}{16y^2}\mod {p^3}.
\endalign$$
\par $(\t{\rm ii})$ If $p\e 3\mod 4$, then
$$\align
&\sum_{k=0}^{(p-1)/2}\f{\b{2k}k^3}{(-8)^k(k+2)} \e\f p9\mod {p^2},
 \\&\sum_{k=0}^{(p-1)/2}\f{\b{2k}k^3}{(-8)^k(k+3)}
 \e \f 1{250}R_1(p)+\f{11}{225}p\mod {p^2},
\\&\sum_{k=0}^{(p-1)/2}\f{\b{2k}k^3}{(-8)^k(k+2)^2}\e
-\f 1{18}R_1(p)+\f p9\mod {p^2},
\\&\sum_{k=0}^{(p-1)/2}\f{\b{2k}k^3}{(-8)^k(k+3)^2}\e
\f{47}{5625}R_1(p)+\f {43}{1125}p\mod {p^2},
\\&\sum_{k=0}^{(p-1)/2}\f{\b{2k}k^3}{(-8)^k(k+1)^3}\e
-1+12R_1(p)+6p\mod {p^2}.
\endalign$$
\endpro

\pro{Conjecture 13.4} Let $p$ be a prime with $p>5$.
\par $(\t{\rm i})$ If $p\e 1,3\mod 8$ and so $p=x^2+2y^2$, then
$$\align
&(-1)^{\f{p-1}2}\sum_{k=0}^{(p-1)/2}\f{\b{2k}k^3}{(-64)^k(k+2)}
 \e  2x^2-\f 89p-\f{2p^2}{9x^2}\mod {p^3},
 \\&(-1)^{\f{p-1}2}\sum_{k=0}^{(p-1)/2}\f{\b{2k}k^3}{(-64)^k(k+3)}
 \e \f{198}{125}x^2-\f {1036}{1125}p-\f{52p^2}{1125x^2}\mod
{p^3},
\\&(-1)^{\f{p-1}2}\sum_{k=0}^{(p-1)/2}\f{\b{2k}k^3}{(-64)^k(k+1)^2}
\e -8y^2-\f{p^2}{y^2}\mod {p^3},
\\&(-1)^{\f{p-1}2}\sum_{k=0}^{(p-1)/2}\f{\b{2k}k^3}{(-64)^k(k+2)^2}
\e -\f{40}9y^2+\f{32}{27}p+\f{p^2}{3y^2}\mod {p^3},
\\&(-1)^{\f{p-1}2}\sum_{k=0}^{(p-1)/2}\f{\b{2k}k^3}{(-64)^k(k+3)^2} \e
\f{56}{5625}y^2-\f{1312}{16875}p-\f{81p^2}{625y^2}\mod {p^3},
\\&(-1)^{\f{p-1}2}\sum_{k=0}^{(p-1)/2}\f{\b{2k}k^3}{(-64)^k(k+1)^3}
\e -8(-1)^{\f{p-1}2}-32y^2+8p-\f{4p^2}{y^2}\mod {p^3}.
\endalign$$
\par $(\t{\rm ii})$ If $p\e 5,7\mod 8$, then
$$\align
&(-1)^{\f{p-1}2}\sum_{k=0}^{(p-1)/2}\f{\b{2k}k^3}{(-64)^k(k+2)} \e\f
7{54}R_2(p)+\f p9\mod {p^2},
\\&(-1)^{\f{p-1}2}\sum_{k=0}^{(p-1)/2}\f{\b{2k}k^3}{(-64)^k(k+3)}
\e -\f {19}{270}R_2(p)-\f {29}{225}p\mod {p^2},
\\&(-1)^{\f{p-1}2}\sum_{k=0}^{(p-1)/2}\f{\b{2k}k^3}{(-64)^k(k+2)^2}
\e\f{19}{27}R_2(p)+\f {2}{27}p\mod {p^2},
\\&(-1)^{\f{p-1}2}\sum_{k=0}^{(p-1)/2}\f{\b{2k}k^3}{(-64)^k(k+3)^2}
\e -\f{233}{675}R_2(p)-\f {254}{3375}p\mod {p^2},
\\&(-1)^{\f{p-1}2}\sum_{k=0}^{(p-1)/2}\f{\b{2k}k^3}{(-64)^k(k+1)^3}
\e -8(-1)^{\f{p-1}2}-12R_2(p)\mod {p^2}.
\endalign$$
\endpro

\pro{Conjecture 13.5} Let $p$ be a prime with $p>5$.
\par $(\t{\rm i})$ If $p\e 1\mod 4$ and so $p=x^2+4y^2$, then
$$\align
&(-1)^{\f {p-1}4}\sum_{k=0}^{(p-1)/2}\f{\b{2k}k^3}{(-512)^k(k+2)}
 \e -\f{152}{27}x^2+\f{190}{27}p-\f{22p^2}{27x^2}\mod {p^3},
\\&(-1)^{\f {p-1}4}\sum_{k=0}^{(p-1)/2}\f{\b{2k}k^3}{(-512)^k(k+3)}
 \e \f{4504}{135}x^2-\f{23074}{675}p+\f{9104p^2}{3375x^2}\mod {p^3},
\\&(-1)^{\f{p-1}4}\sum_{k=0}^{(p-1)/2}\f{\b{2k}k^3}{(-512)^k(k+1)^2}\e
64y^2-8p-\f{p^2}{y^2}\mod {p^3},
\\&(-1)^{\f{p-1}4}\sum_{k=0}^{(p-1)/2}\f{\b{2k}k^3}{(-512)^k(k+2)^2}\e
-\f{6464}{27}y^2+\f{664}{27}p+\f{65p^2}{27y^2}\mod {p^3},
\\&(-1)^{\f{p-1}4}\sum_{k=0}^{(p-1)/2}\f{\b{2k}k^3}{(-512)^k(k+1)^3}\e
-64(-1)^{\f{p-1}4}-192y^2+72p-\f{3p^2}{y^2}\mod {p^3}.
\endalign$$
\par $(\t{\rm ii})$ If $p\e 3\mod 4$, then
$$\align &(-1)^{\f {p-3}4}\sum_{k=0}^{(p-1)/2}\f{\b{2k}k^3}{(-512)^k(k+2)}
 \e \f {28}3R_1(p)+\f{38}9p\mod {p^2},
\\&(-1)^{\f {p-3}4}\sum_{k=0}^{(p-1)/2}\f{\b{2k}k^3}{(-512)^k(k+3)}
 \e -\f {14348}{375}R_1(p)-\f{3938}{225}p\mod {p^2},
\\&(-1)^{\f {p-3}4}\sum_{k=0}^{(p-1)/2}\f{\b{2k}k^3}{(-512)^k(k+2)^2}
 \e \f{424}9R_1(p)-\f{16}3p\mod {p^2},
\\&(-1)^{\f {p-3}4}\sum_{k=0}^{(p-1)/2}\f{\b{2k}k^3}{(-512)^k(k+1)^3}
 \e -64(-1)^{\f{p-3}4}-96R_1(p)+24p\mod {p^2}.
 \endalign$$
 \endpro

 \pro{Conjecture 13.6} Let $p$ be a prime with $p>5$.
\par $(\t{\rm i})$ If $p\e 1\mod 3$ and so $p=x^2+3y^2$, then
$$\align
&(-1)^{\f{p-1}2}\sum_{k=0}^{(p-1)/2}\f{\b{2k}k^3}{256^k(k+2)}
 \e -\f{8}{27}x^2+\f {10}9p+\f{2p^2}{9x^2}\mod {p^3},
 \\&(-1)^{\f{p-1}2}\sum_{k=0}^{(p-1)/2}\f{\b{2k}k^3}{256^k(k+3)}
 \e  -\f{12344}{3375}x^2+\f {4354}{1125}p+\f{728p^2}{1125x^2}\mod
 {p^3},
 \\&(-1)^{\f{p-1}2}\sum_{k=0}^{(p-1)/2}\f{\b{2k}k^3}{256^k(k+1)^2}\e
-48y^2+8p\mod {p^3},
\\&(-1)^{\f{p-1}2}\sum_{k=0}^{(p-1)/2}\f{\b{2k}k^3}{256^k(k+2)^2}\e
-48y^2+\f{200}{27}p\mod {p^3},
\\&(-1)^{\f{p-1}2}\sum_{k=0}^{(p-1)/2}\f{\b{2k}k^3}{256^k(k+1)^3}\e
32(-1)^{\f{p-1}2}+96y^2-24p-\f{2p^2}{y^2} \mod {p^3}.\endalign$$
\par $(\t{\rm ii})$ If $p\e 2\mod 3$, then
$$\align&(-1)^{\f{p-1}2}\sum_{k=0}^{(p-1)/2}\f{\b{2k}k^3}{256^k(k+2)}
\e \f {176}{27}R_3(p)+\f {26}{27}p\mod {p^2},
\\&(-1)^{\f{p-1}2}\sum_{k=0}^{(p-1)/2}\f{\b{2k}k^3}{256^k(k+3)}
\e \f {1808}{135}R_3(p)+\f {1378}{675}p\mod {p^2},
\\&(-1)^{\f{p-1}2}\sum_{k=0}^{(p-1)/2}\f{\b{2k}k^3}{256^k(k+2)^2}
\e 32R_3(p)-\f{16}{27}p\mod {p^2},
\\&(-1)^{\f{p-1}2}\sum_{k=0}^{(p-1)/2}\f{\b{2k}k^3}{256^k(k+3)^2}
\e \f{13856}{225}R_3(p)-\f{1184}{675}p\mod {p^2},
\\&(-1)^{\f{p-1}2}\sum_{k=0}^{(p-1)/2}\f{\b{2k}k^3}{256^k(k+1)^3}
\e 32(-1)^{\f{p-1}2}+128R_3(p)+8p\mod{p^2}.
\endalign$$
\endpro

\pro{Conjecture 13.7} Let $p$ be a prime with $p>7$.
\par $(\t{\rm i})$ If $p\e 1,2,4\mod 7$ and so $p=x^2+7y^2$, then
$$\align
&(-1)^{\f{p-1}2}\sum_{k=0}^{(p-1)/2}\f{\b{2k}k^3}{4096^k(k+2)}
 \e \f{52616}{27}y^2-\f{34}{27}p-\f{20p^2}{27y^2}\mod {p^3},
 \\&(-1)^{\f{p-1}2}\sum_{k=0}^{(p-1)/2}\f{\b{2k}k^3}{4096^k(k+3)}
 \e \f{217125848}{3375}y^2-\f{250882}{3375}p
 -\f{83972p^2}{3375y^2}\mod {p^3},
 \\&(-1)^{\f{p-1}2}\sum_{k=0}^{(p-1)/2}\f{\b{2k}k^3}{4096^k(k+1)^2}\e
-1136y^2+64p+\f{2p^2}{y^2}\mod {p^3},
\\&(-1)^{\f{p-1}2}\sum_{k=0}^{(p-1)/2}\f{\b{2k}k^3}{4096^k(k+1)^3}\e
512(-1)^{\f{p-1}2}+6432y^2-648p-\f{6p^2}{y^2}\mod {p^3}.
\endalign$$
\par $(\t{\rm ii})$ If $p\e 3,5,6\mod 7$, then
$$\align
&(-1)^{\f{p-1}2}\sum_{k=0}^{(p-1)/2}\f{\b{2k}k^3}{4096^k(k+2)}
 \e -1216R_7(p)-\f{11266}9p\mod {p^2},
 \\&(-1)^{\f{p-1}2}\sum_{k=0}^{(p-1)/2}\f{\b{2k}k^3}{4096^k(k+3)}
\e -\f{199232}5R_7(p)-\f{1845802}{45}p\mod {p^2},
\\&(-1)^{\f{p-1}2}\sum_{k=0}^{(p-1)/2}\f{\b{2k}k^3}{4096^k(k+1)^3}
\e 512(-1)^{\f{p-1}2}-1536R_7(p)-1944p\mod {p^2}.
\endalign$$
\endpro

\pro{Conjecture 13.8} Let $p>7$ be a prime. Then
$$\align
&\sum_{k=0}^{p-3}\f{\b{2k}k^2\b{3k}k}{(-192)^k(k+2)} \\&\e \cases
 \f 25x^2-\f
7{15}p-\f{23p^2}{20x^2}\mod {p^3}&\t{if $3\mid p-1$ and
$4p=x^2+27y^2$,}
\\-(2p+1)\b{[2p/3]}{[p/3]}^2-\f p3\mod {p^2}&\t{if $3\mid p-2$,}
\endcases
\\&\sum_{k=0}^{p-4}\f{\b{2k}k^2\b{3k}k}{(-192)^k(k+3)}
\\&\e\cases \f {99}{200}x^2-\f {103}{75}p+\f{141p^2}{700x^2}\mod {p^3}
& \t{if $3\mid p-1$ and $4p=x^2+27y^2$,}
\\\f{13}{14}(2p+1)\b{[2p/3]}{[p/3]}^2+\f {23}{60}p\mod {p^2}
& \t{if $3\mid p-2$,}
\endcases
\\&\sum_{k=0}^{p-2}\f{\b{2k}k^2\b{3k}k}{(-192)^k(k+1)^3}
\\&\e\cases -16+\f
{51}8x^2-9p+\f{147p^2}{4x^2}\mod {p^3}&\t{if $3\mid p-1$ and
$4p=x^2+27y^2$,}\\-16+\f{115}2(2p+1)\b{[2p/3]}{[p/3]}^2-\f
{15}{4}p\mod {p^2} & \t{if $3\mid p-2$}\endcases
\endalign$$
and
$$\sum_{k=0}^{p-2}\f{\b{2k}k^2\b{3k}k}{(-192)^k(k+1)^2} \e \f
14x^2-2p+\f{19p^2}{2x^2}\mod {p^3}\  \t{for}\  3\mid p-1\ \t{and}\
4p=x^2+27y^2.$$
\endpro

\pro{Conjecture 13.9} Let $p>11$ be a prime. Then
$$\align&\sum_{k=0}^{p-3}\f{\b{2k}k^2\b{4k}{2k}}{(-144)^k(k+2)}
\\&\q\e\cases -\f{32}5y^2+\f{16}{15}p+\f{19p^2}{315y^2}\mod {p^3} &\t{if
$p=x^2+3y^2\e 1\mod 3$,}
\\-\f 4{21}R_3(p)\mod {p^2}&\t{if $p\e 2\mod 3$,}
\endcases
\\&\sum_{k=0}^{p-4}\f{\b{2k}k^2\b{4k}{2k}}{(-144)^k(k+3)}
\\&\q\e\cases -\f{3584}{825}y^2+\f{1682}{2475}p+\f{431p^2}{17325y^2}\mod
{p^3} &\t{if $p=x^2+3y^2\e 1\mod 3$,}
\\\f 4{63}R_3(p)+\f 2{45}p\mod {p^2}&\t{if $p\e 2\mod 3$,}
\endcases
\\&\sum_{k=0}^{p-2}\f{\b{2k}k^2\b{4k}{2k}}{(-144)^k(k+1)^3}
\\&\q\e\cases -6-\f{376}{9}y^2+\f{56}{9}p-\f{211p^2}{54y^2}\mod
{p^3} &\t{if $p=x^2+3y^2\e 1\mod 3$,}
\\-6+\f {1360}{27}R_3(p)+\f {20}{27}p\mod {p^2}&\t{if $p\e 2\mod 3$}
\endcases
\endalign$$
and
$$\sum_{k=0}^{p-2}\f{\b{2k}k^2\b{4k}{2k}}{(-144)^k(k+1)^2}
\e -\f{40}3y^2+\f 23p-\f{13p^2}{18y^2}\mod {p^3}\ \t{for}\
p=x^2+3y^2 \e 1\mod 3.$$
\endpro
\par{\bf Remark 13.1} There are many conjectures similar to
Conjectures 13.1-13.9. One may consult [S14, Conjectures 5.1-5.46]
and make similar conjectures by Maple.

\end{document}